\documentclass[12pt]{article}
\usepackage{mathrsfs}
\usepackage{amsmath}
\usepackage{float}
\usepackage{cite}
\usepackage{esint}
\usepackage[all]{xy}
\usepackage{amssymb}
\usepackage{amsfonts}
\usepackage{amsthm}
\usepackage{color}
\usepackage{graphicx}
\usepackage{hyperref}
\usepackage{bm}
\usepackage{indentfirst}
\usepackage{geometry}
\usepackage{tikz-cd}
\theoremstyle{plain}%plain,definition,remark. need amsthm package
\newtheorem{thm}{Theorem}[section]

\newtheorem{defn}[thm]{Definition}
\newtheorem{prop}[thm]{Proposition}

\newtheorem{lem}[thm]{Lemma}

\newtheorem{them}[thm]{Theorem}
\theoremstyle{remark}
\newtheorem{ex}[thm]{Example}
\newtheorem{rmk}[thm]{Remark}

\newcommand{\Z}{\mathbb{Z}}
\newcommand{\Q}{\mathbb{Q}}
\newcommand{\R}{\mathbb{R}}

\newcommand{\C}{\mathbb{C}}

\newcommand{\dd}{\mathrm{d}}

\newcommand{\im}{\mathrm{Im}}


\numberwithin{equation}{section}

\title{Exotic diffeomorphisms and the filling-in problem for complex surfaces}
\author{Qizheng You and Zuyi Zhang}
\date{ }

\begin{document}

\maketitle
\begin{abstract}
     We construct explicit exotic diffeomorphisms on infinitely many simply connected projective complete-intersection surfaces of general type, providing new counterexamples to Donaldson's question on symplectic Torelli groups. We further construct projective degenerations of complex surfaces that admit topological fillings but no smooth fillings, which is a refinement of Friedman's filling-in problem. The key ingredient is the spin number defined in Konno-Lin-Mukherjee-Munoz-Ech\~aniz \cite{konno2025constraints}. The spin number of Brieskorn–Pham monodromies can be expressed explicitly in terms of their spectral numbers, revealing a connection between the spectral genus of surface singularities and smooth four-manifold topology.
\end{abstract}

\section{Introduction}

In \cite{konno2025constraints}, the question asked by Donaldson on symplectic Torelli groups was answered. To be specific, the question is about whether the symplectic Torelli group 
\[
I(X,\omega):=\ker(\pi_0\mathrm{Symp}(X,\omega)\rightarrow\mathrm{Aut}H^*(X;\Z)).
\]
is generated by the squared Dehn–Seidel twists (Equation \ref{equ:dehntwist}) on Lagrangian spheres. Konno-Lin-Mukherjee-Munoz-Ech\~aniz \cite{konno2025constraints} constructed exotic diffeormorphisms on elliptic surfaces. Exotic diffeomorphisms are defined as diffeomorphisms that are topologically isotopic to the identity but not smoothly. These exotic diffeomorphisms are defined by compositions of Dehn-Seidel twists and are not generated by the squared Dehn–Seidel twists, so they provides a negative answer to Donaldson's question. As a result, these diffeomorphism are new elements in $I(X,\omega)-K(X,\omega)$ other than squared Dehn-Seidel twists \cite{seidel2008lectures}, where
\[
K(X,\omega):=\ker(\pi_0\mathrm{Symp}(X,\omega)\rightarrow\pi_0\mathrm{Diff}(X)).
\]

In this paper, we prove the following generalized result:
\begin{them}\label{thm:mainexot}
    There are infinitely many closed simply connected complex surfaces $X$ of general type admitting an exotic diffeomorphism $\varphi$. The space $X$ is a complete intersection of $\mathbb{CP}^4$, equipped with the restricted Fubini-Study form $\omega_{\mathrm{FS}}$ from $\mathbb{CP}^4$. The diffeomorphism $\varphi\in I(X,\omega_{\mathrm{FS}})-K(X,\omega_{\mathrm{FS}})$ is not generated by the squared Dehn–Seidel twists.
\end{them}

The question whether $I(X,\omega_{\mathrm{FS}})=K(X,\omega_{\mathrm{FS}})$ was first proposed by Donaldson. Seidel gave a negative answer to this question using Dehn-Seidel twists and \cite{konno2025constraints} showed that the these twists are not the only elements in $I(X,\omega_{\mathrm{FS}})-K(X,\omega_{\mathrm{FS}})$ if $X$ is an elliptic surface. We enlarge the class of such complex surfaces to surfaces of general type.

The existence of exotic diffeomorphisms on projective complete intersection surfaces are already shown in \cite{baraglia2026irreducible}. However, the exotic diffeomorphisms are not explicitly in \cite{baraglia2026irreducible}. The exotic diffeomorphisms in Theorem \ref{thm:mainexot} are constructed in a concrete way.\\

The criterion to determine the exoticness of diffeomorphisms is from \cite[Theorem B]{konno2025constraints}, presented in Theorem \ref{thm:kllmme}. This theorem is based on the family of Seiberg-Witten theory. One of the difficulty to apply the theorem is to calculate the spin number defined in \cite{konno2025constraints}. The spin number, whose definition can be found in the second half of Section \ref{sec:pre}, is an element of $\pi_1\mathrm{SO}(b^+)$ induced from compositions of Dehn-Seidel twists. Konno-Lin-Mukherjee-Munoz-Ech\~aniz developed a computer program to do the calculation for specific examples. But in this paper, we give an explicit formula:

\begin{them}\label{thm:mainspinnum}
    Given a Brieskorn-Pham surface
    \[
    M_\varepsilon=\{(z_0,z_1,z_2)|\,z_0^{a_0}+z_1^{a_1}+z_2^{a_2}=\varepsilon,\,|z_0|^2+|z_1|^2+|z_2|^2\le\delta\},\quad \varepsilon\in \C,
    \]
    let $S_1,\ldots,S_\mu$ be its vanishing cycles with $\mu:=\dim H^2(M_\varepsilon)$. Then the spin number $[\Delta(S_1,\ldots,S_\mu)^h]$ is expressed as
    \begin{equation}\label{equ:spinformu}
        [\Delta(S_1,\ldots,S_\mu)^h]\equiv h\Sigma_{(i,j,k)\in\mathcal{I}}(\frac{i}{a_0}+\frac{j}{a_1}+\frac{k}{a_2}-1)\pmod2,
    \end{equation}
    where $h=\mathrm{lcm}(a_0,a_1,a_2)$ and 
    \[
    \mathcal{I}:=\{(i,j,k)|\, 1\le i\le a_0, 1\le j\le a_1,1\le k\le a_2, \frac{i}{a_0}+\frac{j}{a_1}+\frac{k}{a_2}<1\}.
    \]
\end{them}

\begin{rmk}
    This theorem is true for all weighted homogeneous complex surfaces using the same method.
\end{rmk}

The right hand side of Equation \ref{equ:spinformu} equals to $-h\Tilde{p_g}$, where the spectral genus $\Tilde{p_g}$ is defined in \cite{eriksson2026spectral} and is considered as a secondary invariant to the geometric genus $p_g$. As a result, the above theorem gives a smooth topological interpretation of the parity of the normalized spectral genus.

Theorem \ref{thm:mainspinnum} follows directly from Proposition \ref{prop:main}. The proof is based on the known results in \cite{kollar2013singularities}, \cite{Steenbrink1977}, and \cite{movasati2007mixed} from algebraic geometry and complex geometry.\\

The spin number turns out to be an obstruction of the filling problems, proposed by Friedman \cite{friedman1983degenerating,friedman1985some,friedman1986simultaneous}. Friedman asked if a holomorphic fibration over a disc with a singular fiber at the center has a holomorphic filling up to base changes (Definition \ref{def:basechange}).
\begin{defn}
    Let $\mathcal{X}\rightarrow D$ be a holomorphic fibration over the disc that has only one singular fiber at the center. If there is a holomorphic fiber bundle $\mathcal{Y}\rightarrow D$ and an isomorphism 
    \[
    \mathcal{Y}|_{D^*}\cong\mathcal{X}|_{D^*},
    \]
    then $\mathcal{X}$ is said to have a holomorphic filling, where $D^*$ is $D$ with center removed.
\end{defn}
\noindent One can define the parallel concepts of fillings with the term ``holomorphic" replaced by ``topological" or ``smooth." Some work related to the holomorphic fillings is done by Morgan \cite{morgan1983topological}. Friedman showed that in most cases, there are no holomorphic fillings. It is interesting to ask whether such fibrations has a smooth filling. In Section \ref{sec:obsfill}, we prove that the spin number is an obstruction to smooth fillings. In fact, there are infinitely many fibrations that admit topological fillings but no smooth fillings.

\begin{them}\label{them:smoofill}
    There are infinitely many fibrations $\mathcal{X}$ over the disc such that
    \begin{itemize}
        \item the total space is smooth,
        \item there is only one singular fiber at the center and the singular fiber has a single singular point,
        \item ${\mathcal{X}^{(h)}}^*$, the $h$-order base change of $\mathcal{X}$, admits a topological filling but no smooth filling.
    \end{itemize}
\end{them}

The organization of the paper is described below. Section \ref{sec:pre} introduces necessary backgrounds including Dehn-Seidel twists and the weighted path. Section \ref{sec:GMconn} provides the definition of Gauss-Manin connection. The eigenvectors of this connection are calculated and this is indispensable for deriving the formula in Theorem \ref{thm:mainspinnum}. Section \ref{sec:spinf} contains the proof of Theorem \ref{thm:mainspinnum}. Section \ref{sec:comemb} constructs embeddings of Brieskorn-Pham surfaces into projective complete intersection surfaces and proves the exotic part of Theorem \ref{thm:mainexot}. Section \ref{sec:sympemb} isotopes the embeddings in the previous section into symplectic embeddings  and complete the proof of Theorem \ref{thm:mainexot}.

\noindent\textbf{Acknowledgments:} The authors would like to thank Jianfeng Lin for helpful discussions and  suggestions.\\

\noindent\textbf{AI Disclaimer:} ChatGPT was used for exploration purposes and grammar checking. The paper is full written by human. The authors take the full responsibility for the correctness.

\section{Preliminary}\label{sec:pre}

\begin{defn}[Brieskorn-Pham surfaces]\label{def:bpsurface}
    A Brieskorn-Pham surface is the resolution of the zero-set of the polynomial $f(z_0,z_1,z_2)=z_0^{a_0}+z_1^{a_1}+z_2^{a_2}$ in a small ball of $\C^3$, for $(z_0,z_1,z_2)\in\C^3$ and $a_0,a_1,a_2\ge2$ that are coprime to each other. To be specific, a Brieskorn-Pham surface is the set 
    \[
    \{(z_0,z_1,z_2)|\,z_0^{a_0}+z_1^{a_1}+z_2^{a_2}=\varepsilon,\,|z_0|^2+|z_1|^2+|z_2|^2\le\delta\},\quad \varepsilon\in \C,
    \]
    where $\delta$ is a small positive number. All surfaces are diffeomorphic for generic small $\varepsilon$. 
\end{defn}

\begin{rmk}
    The additional assumption that $a_0$, $a_1$, and $a_2$ are coprime to each other is a technical condition that guarantees the boundary of Brieskorn-Pham surfaces to be $\Z$-homology spheres. This is used in the last step of proving the main theorem. In fact, it is enough to assume that the boundary of Brieskorn-Pham surfaces are $\Q$-homology spheres. The equivalence condition of this requirement can be found in \cite[Proposition 2.2]{martin2023normal}.
\end{rmk}

\textbf{The next goal is to provide the necessary background for the topology of Brieskorn-Pham surfaces.} We first introduce the vanishing cycle, which are the generators of the second homology of the Brieskorn-Pham surfaces. Throughout the paper, $\mu$ is denoted as the dimension of the second homology of the Brieskorn-Pham surface. In fact, 
\[
\mu=(a_0-1)(a_1-1)(a_2-1),
\]
with the notation from Definition \ref{def:bpsurface} (c.f. \cite{zoladek2006monodromy}).\\

Let $a,b,c\in\C$ be generic parameters such that 
\begin{equation}\label{equ:perbf}
    \Tilde{f}:=f+az_0+dz_1+ez_2:B^3_\delta(0)\rightarrow D_\epsilon(0)\subset\C
\end{equation}
has only nondegenerate critical points in $D_\epsilon(0)$, where $\epsilon$ and $\delta$ are two small positive numbers. Such a function $\Tilde{f}$ is called a Morsification of $f$. Denote the critical values as $s_1,\ldots,s_\mu$ and the unique critical point as $p_i\in\Tilde{f}^{-1}(s_i)\cap B^3_\delta(0)$ with $i=1,\ldots,\mu$. Based on these, one can construct a \textbf{distinguished basis of vanishing path $\gamma_1,\ldots,\gamma_\mu$} for $\Tilde{f}$ such that
\begin{itemize}
    \item $\gamma_i(0)=s_0$, $\gamma_i(1)=s_i$, where $s_0$ is a fixed based point on the boundary of the disc $D_\epsilon(0)$ and $i=1,\ldots,\mu$,
    \item two different paths meet only at $s_0$,
    \item the derivatives $\gamma'_1(0),\ldots,\gamma'_\mu(0)$ are are pairwise distinct, and the ordering of the paths $\gamma_1,\ldots,\gamma_\mu$ is by clockwise outgoing order from $s_0$.
\end{itemize}
With suitably chosen local coordinates $(z_0,z_1,z_2)$ centered around a critical point $p_i$, $\Tilde{f}$ can be written of the form
\[
\Tilde{f}=s_i+z_0^2+z_1^2+z_2^2.
\]
For sufficiently small $t>0$, \textbf{the vanishing sphere of $p_i$} is defined as the sphere
\begin{equation}\label{equ:vansphere}
    S_i(t):=\sqrt{\gamma_i(t)-s_i}S^2,
\end{equation}
where $S^2=\{(z_0,z_1,z_2)\in\C^3|\,\im z_j=0, \ \text{for} \ j=0,1,2, \ \text{and} \ \sum_{j=0}^2 z_j^2=1\}$ is the unit sphere. Notice that $S_i(t)\subset\Tilde{f}^{-1}(\gamma_i(t))\cap B^3_\delta(0)$. The \textbf{vanishing cycle} is the homology class $[S_i]$ of the manifold $\Tilde{f}^{-1}(s_0)\cap B^3_\delta(0)$ obtained by transporting $S_i(t)$ along the path $\gamma_i$ from $\gamma_i(t)$ to $\gamma_i(0)=s_0$. More details can be found in \cite{ebeling2020distinguished}. In fact, the vanishing cycles $[S_i]$ generate the second homology of the Brieskorn-Pham surface and $[S_i]$ can be represented by an embedded sphere, which is still denoted as $S_i$.

\begin{prop}[{\cite[Proposition 1, Theorem 3]{ebeling2020distinguished}}]\label{prop:deltloop}
    Let $M:=\Tilde{f}^{-1}(s_0)\cap B^3_\delta(0)$ be a Brieskorn-Pham surface, where $\Tilde{f}$ is the Morsification of $f$ $($Equation \ref{equ:perbf}$)$. Then the vanishing cycles $[S_i]$ defined above, for $i=1,\ldots,\mu$, generate $H_2(M)$. Moreover, the self-intersection of $S_i$ is $-2$ for all $i$.
\end{prop}

\begin{rmk}
    In \cite{ebeling2020distinguished}, the above proposition is more general in the sense that the Milnor fiber is of complex $n$-dimensional. In this paper, we are restricted to the case where $n=2$.
\end{rmk}

\textbf{Then we construct a self-diffeomorphism that is proved to be exotic later.} With an embedded sphere $S\hookrightarrow M$ that has self-intersection number being $-2$, it is possible to define the \textit{Dehn-Seidel twist} $\tau_S$, c.f. \cite{seidel2008lectures}, \cite{scorpan2005wild}, or \cite{konno2025constraints}. The induced action on $H^2(M)$ is given by the \textit{Picard–Lefschetz reflection}:
\begin{equation}\label{equ:dehntwist}
    \tau_S^*(\alpha)=\alpha+\alpha([S])\mathrm{PD}([S]), \quad\alpha\in H^2(M).
\end{equation}
Here $[S]$ refers to the element in $H_2(M)$ represented by the embedded sphere $S$. The action $\tau_S^*$ is a reflection of $H^2(M)$ along the hyperplane perpendicular to $\mathrm{PD}(S)$, and hence we have $(\tau_S^*)^2=\mathrm{id}_{H_2(M)}$. 

Suppose the vanishing cycles of a Brieskorn-Pham surface $M$ are $[S_1],\ldots,[S_\mu]$, then the composition of the Dehn-Seidel twists $\tau_{S_1}\circ\ldots\circ\tau_{S_\mu}$ is a self-diffeomorphism of $M$. The induced map on $H^2(M)$ has a finite order that is a divisor of $h$. If $M$ can be embedded into a closed 4-manifold $X$ such that the extension of $(\tau_{S_1}\circ\ldots\circ\tau_{S_\mu})^h$ acts trivially on $H^2(X)$, then $(\tau_{S_1}\circ\ldots\circ\tau_{S_\mu})^h$ is topologically isotopic to the identity on $X$ by \cite{quinn1986isotopy}. We want to show that the diffeomorphism $(\tau_{S_1}\circ\ldots\circ\tau_{S_\mu})^h$ is not smoothly isotopic to the identity.

To achieve this goal, two paths in $\mathrm{SO}(b^+)$ are needed, where $b^+$ is the dimension of the positive definite subspace of $H^2(M)$. These two paths are proved to be homotopic when $b^+>2$ in the next section.

\begin{figure}[H] 
\centering 
\includegraphics[width=0.5\textwidth]{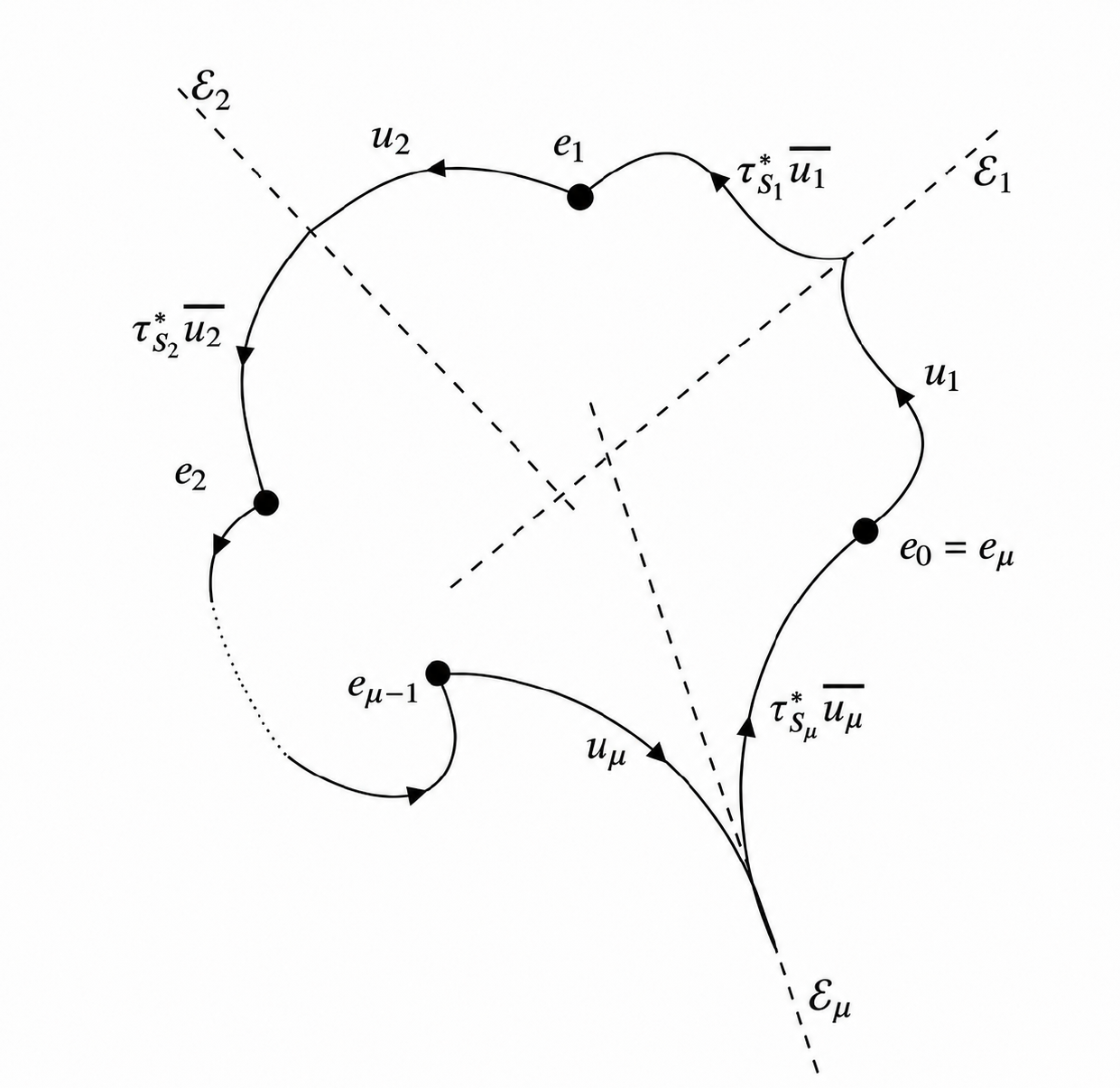}
\caption{This is a picture illustrating the loop of vanishing cycle from \cite{konno2025constraints} with $n=\mu$.}
\label{fig:loopsket}
\end{figure}

\textbf{The first path} is given by \cite{konno2025constraints} and we describe its construction below. Let $\mathcal{E}$ be the space of linear embeddings $e:\R^{b^+}\hookrightarrow H^2(M)$ whose image $\im(e)$ is a positive subspace with respect to the intersection form of $H^2(M)$. A given embedding $e_0\in\mathcal{E}$ yields a homotopy equivanlence $\mathrm{SO}(b^+)\simeq\mathcal{E}$ by reparametrizaiton of $e_0$. Let $e_i=\tau_{S_i}^*\circ\ldots\circ\tau_{S_1}^*\circ e_0$ for $i\leq\mu$, and $u_i$ be a path in $\mathcal{E}$ from $e_{i-1}$ to the subspace $\mathcal{E}_i\subset\mathcal{E}$ consisting of embeddings whose image is orthogonal to $S_i$. Then the concatenation of the $2\mu$ paths:
\[
u_1,\tau_{S_1}^*\circ\overline{u_1},u_2,\tau_{S_2}^*\circ\overline{u_2},\ldots,u_\mu,\tau_{S_\mu}^*\circ\overline{u_\mu}
\]
is an element of $\pi_1\mathrm{SO}(b^+)$, where $\overline{u_i}$ stands for the reversed path of $u_i$. This path is denoted by $\Delta(S_1,\ldots,S_\mu)$ and is called \textbf{the path of vanishing cycles} (Figure \ref{fig:loopsket}). In fact, $\Delta(S_1,\ldots,S_\mu)\in\mathrm{SO}(b^+)$ is well-defined up to the choice of $e_0,u_1,\ldots,u_\mu$. Furthermore, concatenating $\Delta(S_1,\ldots,S_\mu)$ with itself for $h$ times, we get a loop $\Delta(S_1,\ldots,S_\mu)^h$ in $\mathrm{SO}(b^+)$. The induced element in $\pi_1\mathrm{SO}(b^+)$ is called the \textit{spin number} in \cite{konno2025constraints}.

\begin{lem}[{\cite[Lemma 2.4, Proposition 2.8]{konno2025constraints}}]
    The element $[\Delta(S_1,\ldots,S_\mu)]\in\mathrm{SO}(b^+)$ is independent of the auxiliary choices of $e_0,u_1,\ldots,u_\mu$. And one representation is given by the concatenating of the following $\mu$ paths:
    \[
    \rho_t^i\circ e_{i-1},\ 0 \leq t \leq 1, \ i=1,\ldots,\mu,
    \]
    where
    \begin{align*}
        \rho_t^i: H^2(M)&\rightarrow H^2(M)\\
        [\alpha]&\mapsto \alpha+t\alpha([S_i])PD([S_i]).
    \end{align*}
\end{lem}

\textbf{The second path} is defined on the fiber of $f$ instead of $\Tilde{f}$. Recall that $M$ is diffeomorphic to ${f}^{-1}(\varepsilon)\cap B^3_\delta(0)$ with $\varepsilon$ being a complex number in the unit circle \cite{milnor1968singular}. For $t \in D_\epsilon(0)$, we have 
\[
f(t^{w_0}z_0,t^{w_1}z_1,t^{w_2}z_2)=t^hf(z_0,z_1,z_2).
\]
Here $h=\mathrm{lcm}(a_0,a_1,a_2)$ is the order of the automorphism $\tau_{S_1}^*\circ\ldots\circ\tau_{S_\mu}^*$ of $H^2(M)$, and we take $w_0,w_1,w_2$ such that $a_0w_0=a_1w_1=a_2w_2=h$. Therefore, along the boundary of the unit disc, the weighted action $(z_0,z_1,z_2)\mapsto(t^{w_0}z_0,t^{w_1}z_1,t^{w_2}z_2)$ induces a diffeomorphism $f^{-1}(1)\rightarrow f^{-1}({t^h})$ for each $t\in S^1$. As $t$ going along the boundary of the unit disc $S^1$, the second cohomology group can be pulled back from the fiber $f^{-1}(t^h)$ to the fiber $f^{-1}(1)$. This induces a path in $\mathrm{SO}(b^+)$ as $t$ goes from $1$ to $\mathrm{e}^{\frac{2\pi \mathrm{i}}{h}}$ along the boundary $S^1$. This path is denoted as $\rho_{wt}$ and is called the \textbf{weighted path}. Furthermore, if $t$ travels along the full circle once, we get a loop in $\mathrm{SO}(b^+)$ and it is denoted as $\rho_{wt}^h$. Detailed description of the Gauss-Manin connection is provided in the next section.

\section{The Gauss-Manin connection}\label{sec:GMconn}

Let $\Tilde{f}_u(z_0,z_1,z_2):=f(z_0,z_1,z_2)+u(az_0+dz_1+ez_2)$ be a path connecting $f$ and $\Tilde{f}$ (Equation \ref{equ:perbf}) for $u\in[0,1]$. For generic $\varepsilon$ in the unit disc and $u$, $\Tilde{f}_u^{-1}(\varepsilon)$ is a smooth Brieskorn-Pham surface. Therefore, by varying $\varepsilon$ for each fixed $u$, we get a family of Brieskorn-Pham surfaces defined on the punctured disc whose punctures are the critical values of $\Tilde{f}_u$. The fiber of this bundle is denoted by $M_{u,\varepsilon}$. Since $M_{u,\varepsilon}$ is diffeomorphic for different $(u,\varepsilon)$ whenever $M_{u,\varepsilon}$ is smooth, after taking the second cohomology of $M_{u,\varepsilon}$, there is a vector bundle $\mathcal{H}$ defined over the punctured disc with fiber $H^2(M_{u,\varepsilon};\C)$ for each fixed $u$.

%The goal of this section is to find a concrete framing of $\mathcal{H}$ using the Poincar\'e residue map. and show that the sections that representing this trivilization can be extended across the punctures. As a consequence, $\mathcal{H}$ can be defined over the whole unit disc. Throughout this section, the numbers $a,b,c$ are assumed to be generic with small enough norms such that all the critical values of $\Tilde{f}_u$ are isolated except $u=0$.\\

The first purpose is to construct a framing of $\mathcal{H}$ using the Poincar\'e residue map \cite{griffiths2014principles}. We review the following $2$-form on each smooth $M_{u,\varepsilon}$ introduced in \cite[Page 147]{griffiths2014principles}. Let
\begin{equation}
    \omega_{\nu,u,\varepsilon}:=(-1)^l\frac{z_0^{i-1}z_1^{j-1}z_2^{k-1}\dd z_0\wedge \dd z_1\wedge \dd z_2}{(\partial\Tilde{f}_u/\partial z_l)\dd z_l}|_{\Tilde{f}_u=\varepsilon},
\end{equation}
where $\nu=(i,j,k)$, $1\le i\le a_0$, $1\le j\le a_1$, $1\le k\le a_2$, and $l=0,1,2$. Note that the $2$-form is independent of the choice of $l$ as long as $\partial\Tilde{f}_u/\partial z_l\ne0$ for any $l=0,1,2$ \cite[Page 147]{griffiths2014principles}. It remains to show that $\{[\omega_{\nu,u,\varepsilon}]\}_\nu$ forms a basis of $H^2(M_{u,\varepsilon};\C)$ for certain $\varepsilon$ and $u$.

\begin{prop}
    For smooth Brieskorn-Pham surfaces $M_{u,\varepsilon}$, the family $\{[\omega_{\nu,u,\varepsilon}]\}_\nu$ is a basis for $H^2(M_{u,\varepsilon};\C)$.
\end{prop}

\begin{proof}
    In \cite{movasati2008moduli}, such proposition is claimed to be correct using the results from \cite{movasati2007mixed} but unproved, so we present a brief proof here. Regarding the notation, we remark here that the polynomial $f$ in \cite{movasati2007mixed} corresponds to $\Tilde{f}_u$ in our paper, so the function $f$ in \cite{movasati2007mixed} is identified with $\Tilde{f}_u$ in our proof. The number $n$ in \cite{movasati2007mixed} is taken to be 2 in our proof.

    In \cite{movasati2007mixed}, Movasati defines the following two $\C[t]$-modules
    \begin{equation}\label{equ:brieskornlat}
        H':=\frac{\dd\Tilde{f}_u\wedge\Omega^n}{\dd\Tilde{f}_u\wedge \dd\Omega^{n-1}},\ \ H'':=\frac{\Omega^{n+1}}{\dd\Tilde{f}_u\wedge \dd\Omega^{n-1}},
    \end{equation}
    where $\Omega^i$ is the the space consisting of polynomial differential $i$-forms in $\C^{n+1}$ for $i=n-1,n,n+1$. The $\C[t]$-module structure is induced by the action $t[\omega]=[\Tilde{f}_u\omega]$, where $[\omega]\in H'$ or $H''$. Using the Atiyah-Hodge-Grothendieck theorem in the introduction of \cite{movasati2007mixed}, Movasati proved that the restriction of $H'$ to each regular fiber $\Tilde{f}_u^{-1}(\varepsilon)$ is isomorphic to the $n$-th de Rham cohomology of $\Tilde{f}_u^{-1}(\varepsilon)$ with complex coefficients. That is to say,
    \[
    H'|_{t=\varepsilon}\cong \frac{H'}{(t-\varepsilon)H'}\cong H^n(M_{u,\varepsilon};\C).
    \]
    It is proved in \cite[Section 3]{movasati2007mixed}, that the localizations of $H'$ and $H''$ at $t-\varepsilon$ are isomorphic through the map
    \begin{align*}
        \rho:H''_{(t-\varepsilon)}&\rightarrow H'_{(t-\varepsilon)}\\
        [\omega]&\mapsto[\frac{\omega}{\dd\Tilde{f}_u}]
    \end{align*}
    if $\varepsilon$ is not a critical value of $\Tilde{f}_u$. Here $\frac{\omega}{d\Tilde{f}_u}$ is defined to be the form $\alpha$ such that $d\Tilde{f}_u\wedge\alpha=\omega$. Because of the following isomorphisms of the restrictions of the localizations at a specific fiber,
    \[
    H''_{(t-\varepsilon)}|_\varepsilon\cong\frac{H''}{(t-\varepsilon)H''},\ \ H'_{(t-\varepsilon)}|_\varepsilon\cong\frac{H'}{(t-\varepsilon)H'},
    \]
    it holds that
    \[
    H''_{(t-\varepsilon)}|_\varepsilon\cong H^n(M_{u,\varepsilon};\C).
    \]
    According to \cite[Proposition 6.2]{movasati2007mixed}, $z_0^{i-1}z_1^{j-1}z_2^{k-1}\dd z_0\wedge \dd z_1\wedge \dd z_2$ is a $\C[t]$-basis of $H''$ for $1\le i\le a_0$, $1\le j\le a_1$, $1\le k\le a_2$. As a result, $\rho$ induces a basis for $H^2(M_{u,\varepsilon};\C)$. That is to say,
    \[
    [\omega_{\nu,u,\varepsilon}]=[\frac{z_0^{i-1}z_1^{j-1}z_2^{k-1}\dd z_0\wedge \dd z_1\wedge \dd z_2}{\dd\Tilde{f}_u}|_{\Tilde{f}_u=\varepsilon}]
    \]
    is the basis of $H^2(M_{u,\varepsilon};\C)$. The above equation holds since by define, $\frac{z_0^{i-1}z_1^{j-1}z_2^{k-1}\dd z_0\wedge \dd z_1\wedge \dd z_2}{\dd\Tilde{f}_u}$ is equal to $\alpha$ such that $d\Tilde{f}_u\wedge=z_0^{i-1}z_1^{j-1}z_2^{k-1}\dd z_0\wedge \dd z_1\wedge \dd z_2$ and $\alpha$ can be taken to be $(-1)^l\frac{z_0^{i-1}z_1^{j-1}z_2^{k-1}\dd z_0\wedge \dd z_1\wedge \dd z_2}{(\partial\Tilde{f}_u/\partial z_l)\dd z_l}$ with direct calculation as long as $\partial\Tilde{f}_u/\partial z_l\ne0$.
\end{proof}

The two paths $\Delta(S_1,\ldots,S_\mu)$ and $\rho_{wt}$ defined at the end of the last section are connected by the Gauss-Manin connection defined as below.

\begin{defn}[\cite{movasati2007mixed}]\label{GSconn}
    Given $u\in[0,1]$, let $\Tilde{f}_u(z_0,z_1,z_2):=f(z_0,z_1,z_2)+u(az_0+bz_1+cz_2)$. The corresponding \textbf{Gauss-Manin connection} is defined as
    \[
    \nabla_{GM}:H'\rightarrow H'',\ \nabla_{GM}([\dd\Tilde{f}_u\wedge\omega])=[\dd\omega],
    \]
    where $H'$ and $H''$ are defined in Equation \ref{equ:brieskornlat}.
\end{defn}

\begin{rmk}
    In \cite[Page 787]{movasati2007mixed}, it is claimed that the map
    \begin{align*}
        \frac{\Omega^n}{\dd \Tilde{f}_u\wedge\Omega^{n-1}+\dd\Omega^{n-1}}&\rightarrow\frac{\dd\Tilde{f}_u\wedge\Omega^n}{\dd\Tilde{f}_u\wedge \dd\Omega^{n-1}}\\
        [\omega]&\mapsto[\dd\Tilde{f}_u\wedge\omega]
    \end{align*}
    is an isomorphism. Moreover, $\frac{\Omega^n}{\dd \Tilde{f}_u\wedge\Omega^{n-1}+\dd\Omega^{n-1}}$ can be identified with the $n$-th de Rham cohomology after restricting to the regular fiber. In fact, this is realized by the map
    \begin{align*}
        \frac{\Omega^n}{\dd \Tilde{f}_u\wedge\Omega^{n-1}+\dd\Omega^{n-1}}&\rightarrow H^2(M_{u,\varepsilon};\C)\\
        [\omega]&\mapsto[\omega|_{M_{u,\varepsilon}}].
    \end{align*}
    For the case $n=2$, the Gauss-Manin connection can act on the second cohomology of the regular fiber through these identifications.
\end{rmk}

We end this section with two examples that illustrate why the Gauss–Manin connection is related to $\Delta(S_1,\ldots,S_\mu)$ and $\rho_{wt}$.

\begin{ex}
    As claimed in \cite[Section 2]{greuel2001applications}, the monodromy of the Gauss-Manin connection along a small meridian around the critical value of $\Tilde{f}_u$ is the Picard-Lefschetz reflection. Details can be found in \cite{brieskorn1970monodromie} \cite{greuel2001applications}.
\end{ex}

\begin{ex}\label{exp:weightedloop}
    Recall that the weighted path is defined by the weighted action $\mathrm{e}^{i\theta}\cdot(z_0,z_1,z_2)=(\mathrm{e}^{iw_0\theta}z_0,\mathrm{e}^{iw_1\theta}z_1,\mathrm{e}^{iw_2\theta}z_2)$. The induced diffeomorphism between the fibers $M_{0,1}$ and $M_{0,\mathrm{e}^{ih\theta}}$ yields an isomorphism between $H^2(M_{0,1};\C)$ and $H^2(M_{0,\mathrm{e}^{ih\theta}};\C)$, where $l=\mathrm{lcm}(a_0,a_1,a_2)$. This isomorphism can be expressed explicitly using the Gauss-Manin connection.

    In \cite[Section 3]{movasati2007mixed}, it is calculated that $\varepsilon^{-\frac{N_\nu}{h}}\omega_{\nu,0,\varepsilon}$ is parallel under the Gauss-Manin connection, where $\varepsilon$ is on the boundary of the unit disc and 
    \[
    N_\nu:=iw_0+jw_1+kw_2-h.
    \]
    In fact, the notations $\nabla$, $A_\beta$, and $\frac{\dd f}{f}\wedge\frac{\eta_\beta}{f-\varepsilon}$ in \cite{movasati2007mixed} correspond to the Gauss-Manin connection, $\frac{N_\nu+h}{h}$, and $\frac{z_0^{i-1}z_1^{j-1}z_2^{k-1}\dd z_0\wedge\dd z_1\wedge z_2}{f-\varepsilon}$, respectively. After applying the Poincare residue map, $\omega_{\nu,0,\varepsilon}$ is in correspondence with $\frac{\eta_\beta}{\varepsilon}|_{M_0,\varepsilon}$ Mosavati's notation. From the equation $\nabla\eta_\beta=\frac{A_\beta}{t}\eta_\beta$ in \cite[Section 3]{movasati2007mixed} and the Leibniz rule, we have (the variable is $\varepsilon$)
    \begin{align*}
        \nabla_{GM}(\varepsilon^{-\frac{N_\nu}{h}}\omega_{\nu,0,\varepsilon})=&-\frac{N_\nu}{h}\varepsilon^{-\frac{N_\nu+h}{h}}\omega_{\nu,0,\varepsilon}+\varepsilon^{-\frac{N_\nu}{h}}\nabla_{GM}\omega_{\nu,0,\varepsilon}\\
        =&-\frac{N_\nu}{h}\varepsilon^{-\frac{N_\nu+h}{h}}\omega_{\nu,0,\varepsilon}+\varepsilon^{-\frac{N_\nu}{h}}\frac{N_\nu}{h\varepsilon}\omega_{\nu,0,\varepsilon}\\
        =&0.
    \end{align*}
    In other words, 
    \[
    P_{0,\varepsilon}^{GM}(\omega_{\nu,0,1})=\varepsilon^{-\frac{N_\nu}{h}}\omega_{\nu,0,\varepsilon},
    \]
    where $P_{0,\varepsilon}^{GM}$ is the Gauss-Manin parallel transport from $H^2(M_{0,1};\C)$ to $H^2(M_{0,\varepsilon};\C)$. As a result,
    \[
    (P_{0,\varepsilon}^{GM})^{-1}(\omega_{\nu,0,\mathrm{e}^{ih\theta}})=\mathrm{e}^{iN_\nu\theta}\omega_{\nu,0,1}.
    \]
    \textbf{This proves that under the $[\omega_{\nu,0,\varepsilon}]$ framing of $H^2(M_{0,\varepsilon};\C)$, the framings along the boundary loop $\mathrm{e}^{ih\theta}$ are of a diagonal form with each diagonal entry being $\mathrm{e}^{iN_\nu\theta}$ along the boundary loop $\mathrm{e}^{ih\theta}$.}
\end{ex}

\section{Sipn number formula}\label{sec:spinf}

The goal of this section is to prove the following proposition and thus proves Theorem \ref{thm:mainspinnum}.

\begin{prop}\label{prop:main}
    Let $M$ be a Brieskorn-Pham surface. Then the two paths $\Delta(S_1,\ldots,S_\mu)$ and $\rho_{wt}$ defined at the end of the last section are homotopic. Moreover, the loops $\Delta(S_1,\ldots,S_\mu)^h$ and $\rho_{wt}^h$ satisfy the following relation:
    \begin{itemize}
        \item $[\Delta(S_1,\ldots,S_\mu)^h]\equiv[\rho_{wt}^h]\equiv \Sigma_{\nu\in\mathcal{I}}{N_\nu}$ $(\mathrm{mod}$ $2)$ in $\pi_1\mathrm{SO}(b^+)$.%$\cong\Z$ for $b^+=2$,
        %\item $[\Delta(S_1,\ldots,S_\mu)]=[\rho_{wt}]$ in $\pi_1\mathrm{SO}(b^+)\cong\Z_2$ for $b^+>2$,
    \end{itemize}
    where $[S_i]$ are the vanishing cycles of $M$ with $i=1,\ldots,\mu$, $\mathcal{I}$ is defined in the lemma below, and $N_\nu:=iw_0+jw_1+kw_2-h$.
\end{prop}

In the last section, we construct a frame for $H^2(M_{u,\varepsilon};\C)$. But what is actually needed is the frame of the maximal positive definite subspace. The next lemma shows that the low-weighted Poincar\'e residue forms gives such a frame.

\begin{lem}\label{lem:posframe}
    Let $\Tilde{f}_u(z_0,z_1,z_2):=f(z_0,z_1,z_2)+u(az_0+bz_1+cz_2)$ and let $M_{u,\varepsilon}$ be $\Tilde{f}^{-1}_u(\varepsilon)$ for a regular value $\varepsilon$ in the unit disc with $u\in[0,1]$. Then $[\mathrm{Re}\omega_{\nu,u,\varepsilon}]$ and $[\mathrm{Im}\omega_{\nu,u,\varepsilon}]$ with
    \[
    \nu\in\mathcal{I}:=\{(i,j,k)|\, 1\le i\le a_0, 1\le j\le a_1,1\le k\le a_2, \frac{i}{a_0}+\frac{j}{a_1}+\frac{k}{a_2}<1\}
    \]
    generates a maximal positive definite subspace of $H^2(M_{u,\varepsilon};\R)$. The forms $\omega_{\nu,u,\varepsilon}$ with $\nu\in\mathcal{I}$ are called the \textbf{low-weighted Poincar\'e residue forms}.
\end{lem}

\begin{proof}
    According to \cite[Theorem 2]{Steenbrink1977}, $b^+=2|\mathcal{I}|$ for Brieskorn-Pham surfaces. The key point is to show that the real and imaginary parts of the low-weighted forms generate a positive definite subspace. The idea is to compactify $M_{u,\varepsilon}$ into a closed K\"ahler manifold $Z_{u,\varepsilon}$ by adding serval divisors and then extend the low-weighted forms to holomorphic two forms $\Omega_{\nu,u,\varepsilon}$ on $Z_{u,\varepsilon}$. The holomorphicity of $\Omega_{\nu,u,\varepsilon}$ on closed K\"ahler manifolds yields the desired property.

    The definition of $Z_{u,\varepsilon}$ follows from \cite{ohta1999simple} \cite{lee2026instanton}. First, $M_{u,\varepsilon}$ can be compactified in the weighted projective space $\mathbb{CP}(w_0,w_1,w_2,1)$ with the formula
    \[
    \overline{M_{u,\varepsilon}}:=\{[z_0:z_1:z_2:z_3]|\ f(z_0,z_1,z_2)+u(az_0z_3^{h-w_0}+bz_1z_3^{h-w_1}+cz_2z_3^{h-w_2})-\varepsilon z_3^h=0\},
    \]
    where $\mathrm{lcm}(a_0,a_1,a_2)=h=a_0w_0=a_1w_1=a_2w_2$. Since $\mathbb{CP}(w_0,w_1,w_2,1)$ is an algebraic projective variety, $\overline{M_{u,\varepsilon}}$ is a complex algebraic variety by adding a divisor at infinity to $M_{u,\varepsilon}$. Then there is a standard way to resolve the singularities of $\overline{M_{u,\varepsilon}}$ through finite times of blow ups by Hironaka \cite{hironaka1964resolution1,hironaka1964resolution} (see also \cite{abramovich2018resolution}) to get a smooth projective surface $\pi:Z_{u,\varepsilon}\rightarrow M_{u,\varepsilon}$. As a result, $Z_{u,\varepsilon}$ is a closed K\"ahler manifold. Moreover, $Z_{u,\varepsilon}-M_{u,\varepsilon}$ is the union of serval smooth divisors with normal crossings according to Hironaka's construction \cite{hironaka1964resolution1,hironaka1964resolution}.% Denote the disjoint union of these divisors as $\hat{D}$.

    The next step is to extend the low-weighted Poincar\'e residue forms to $Z_{u,\varepsilon}$. Without loss of generality, we assume that $\partial\Tilde{f}_u/\partial z_0\ne0$ in some neighborhood $U\subset M_{u,\varepsilon}$. Then it holds in $U$ that
    \[
    \omega_{\nu,u,\varepsilon=}\frac{z_0^{i-1}z_1^{j-1}z_2^{k-1}\dd z_1\wedge \dd z_2}{\partial\Tilde{f}_u/\partial z_0}.
    \]
    Consider the 2-form
    \begin{align*}
            \Omega_{wt}^\nu:=z_0^{i-1}z_1^{j-1}z_2^{k-1}z_3^{r_\nu}(&-w_0z_0\dd z_1\wedge \dd z_2\wedge \dd z_3+w_1z_1\dd z_0\wedge \dd z_2\wedge \dd z_3\\
            &-w_2z_2\dd z_0\wedge \dd z_1\wedge \dd z_3+z_3\dd z_0\wedge \dd z_1\wedge \dd z_2),
    \end{align*}
    where $\nu=(i,j,k)$ with $1\le i\le a_0, 1\le j\le a_1,1\le k\le a_2$ and $r_\nu=h-iw_0-jw_1-kw_2-1$. Since $\frac{i}{a_0}+\frac{j}{a_1}+\frac{k}{a_2}<1$, $r_\nu\ge0$. Then using the Poincar\'e residue map, define
    \[
    \overline{\omega_{\nu,u,\varepsilon}}:=(-1)^l\frac{\Omega_{wt}^\nu}{\partial(F_{u,\varepsilon}(z_0,z_1,z_2,z_3))/\partial z_ldz_l}|_{F_{u,\varepsilon}=0},
    \]
    where $F_{u,\varepsilon}(z_0,z_1,z_2,z_3):=f(z_0,z_1,z_2)+u(az_0z_3^{h-w_0}+dz_1z_3^{h-w_1}+ez_2z_3^{h-w_2})-\varepsilon z_3^h$ and $\overline{\omega_{\nu,u,\varepsilon}}$ is is independent of the choice of $l$ as long as $\partial F_{u,\varepsilon}/\partial z_l\ne0$ for any $l=0,1,2,3$ \cite[Page 147]{griffiths2014principles}. Note that $M_{u,\varepsilon}$ is by setting $z_3=1$ in $\overline{M_{u,\varepsilon}}$, a direct calculation using definitions for $z_3=1$ shows that
    \[
    \overline{\omega_{\nu,u,\varepsilon}}|_{M_{\nu,\varepsilon}}=\omega_{\nu,u,\varepsilon}.
    \]
    Finally, the extension $\Omega_{\nu,u,\varepsilon}$ of $\omega_{\nu,u,\varepsilon}$ to $Z_{u,\varepsilon}$ is defined as the pull-back of $\overline{\omega_{\nu,u,\varepsilon}}$ through the blow up map $\pi:Z_{u,\varepsilon}\rightarrow \overline{M_{u,\varepsilon}}$, i.e.
    \[
    \Omega_{\nu,u,\varepsilon}:=\pi^*\overline{\omega_{\nu,u,\varepsilon}}.
    \]
    {According to the next lemma, $\Omega_{\nu,u,\varepsilon}$ can be extended to a holomorphic 2-form on $Z_{u,\varepsilon}$, still denoted as $\Omega_{\nu,u,\varepsilon}$.} Therefore, it is a fundamental fact for closed K\"ahler 4-manifolds that the space generated by $[\mathrm{Re}\Omega_{\nu,u,\varepsilon}]$ and $[\mathrm{Im}\Omega_{\nu,u,\varepsilon}]$ is positive definite. So, after restricting this space to $M_{u,\varepsilon}$, it is still positive definite by calculating the integral directly. This space is maximal because of the dimension reason.
\end{proof}

Suppose that $\varepsilon$ is a critical value, the 2-form $\omega_{\nu,u,\varepsilon}$ is well defined except at the critical point. Moreover, $\omega_{\nu,u,\varepsilon}$ can be extended over the divisor at the infinity to get $\overline{\omega_{\nu,u,\varepsilon}}$. Therefore, the question that whether $\Omega_{\nu,u,\varepsilon}:=\pi^*\overline{\omega_{\nu,u,\varepsilon}}$ can be extended over the resolved infinity divisor or not is irrelevant to if $\varepsilon$ is a regular value or not.

\begin{lem}\label{lem:extinfty}
    Let $\pi:Z_{u,\varepsilon}\rightarrow \overline{M_{u,\varepsilon}}$ be the log resolution of the compactification of $M_{u,\varepsilon}$ over the divisor at the infinity. Then $\Omega_{\nu,u,\varepsilon}:=\pi^*\overline{\omega_{\nu,u,\varepsilon}}$ can be extended as a holomorphic 2-form across the resolution, where $\varepsilon$ is possibly a critical value.
\end{lem}

\begin{proof}
    The idea is to apply the results from \cite{kollar2013singularities} and \cite{greb2011differential}. The Observation 1.3 and Theorem 1.4 of \cite{greb2011differential} combine to get that for a complex dimension $n$ complex
    quasi-projective variety with klt (Kawamata log terminal) singularities, the pullback of the holomorphic $n$-forms to its log resolutions can be extended holomorphically. Since $M_{u,\varepsilon}$ is projective, only need to show that it is klt.
    
    The criterion for being klt is from \cite[Proposition 5.20]{kollar2013singularities}. This proposition states that if there is a finite morphism $g:V\rightarrow \overline{M_{u,\varepsilon}}$  with $V$ being a smooth complex manifold, and 
    \[
    K_V=g^*K_{\overline{M_{u,\varepsilon}}}
    \]
    between the canonical divisors (i.e. the divisor of a rational 2-form), then $\overline{M_{u,\varepsilon}}$ is klt. {In fact, we only need to construct a finite cover $g:V\rightarrow\overline{M_{u,\varepsilon}}$ around the infinity singularities.}

    It remains to define $V$ and to show that the relation between the canonical divisors. The finite cover is constructed in the coordinates with $z_0=1,z_1=1,$ and $z_2=1$, respectively. Since the these constructions are essentially the same, we only describe the one in $z_0=1$ coordinate $U$. In this coordinate, $\overline{M_{u,\varepsilon}}$ is determined by the equation
    \begin{equation}
        1+z_1^{a_1}+z_2^{a_2}+u(az_3^{h-w_0}+dz_1z_3^{h-w_1}+ez_2z_3^{h-w_2})-\varepsilon z_3^h=0.
    \end{equation}
    In fact, the solution of the above equation in $\C^3$, denoted as $V$, quotient by the finite $\Z_{w_0}$-action
    \[
    \xi\cdot(z_1,z_2,z_3):=(\xi^{w_1}z_1,\xi^{w_2}z_2,\xi z_3),\ \xi\in\Z_{w_0},
    \]
    is the corresponding part of $\overline{M_{u,\varepsilon}}$. Therefore, this space $V$ is the finite covering. It is clear that $V$ is smooth near the infinity divisor ($z_3=0$). {Since $g:V\rightarrow U\cap \overline{M_{u,\varepsilon}}$ is a finite covering, $g$ is a separable finite surjective morphism of normal varieties \cite[Appendix 1, Proposition 1.1]{MustataFiniteGroupQuotients}.} According to \cite[Proposition 2.2]{cascini2020polarized},
    \[
    K_V=g^*K_{\overline{M_{u,\varepsilon}}}|_U+R_g,
    \]
    where $R_g$ is a divisor that vanishes when $g$ is \'etale in codimension 1 (i.e. quasi-\'etale) \cite{cascini2020polarized}. Because the $\Z_{w_0}$-action on $V$ is free except at finitely points, then $g$ is quasi-\'etale. Therefore,
    \[
    K_V=g^*K_{\overline{M_{u,\varepsilon}}}|_U.
    \]
    This argument holds for $z_1=1$ and $z_2=1$ coordinates as well, so $\overline{M_{u,\varepsilon}}$ is klt.
    
    %To verify the relation of the canonical divisors, it is sufficient to show the same relation for the 2-forms defining these divisors. Using the Poincar\'e residue map, in the $z_0=1$ coordinates, we can define a 2-form
    %\[
    %(-1)^l\frac{-w_0dz_1\wedge dz_2\wedge dz_3}{\partial(F_{u,\varepsilon}(1,z_1,z_2,z_3))/\partial z_ldz_l}|_{F_{u,\varepsilon}=0},
    %\]
    %where $F_{u,\varepsilon}(z_0,z_1,z_2,z_3):=f(z_0,z_1,z_2)+u(az_0z_3^{h-w_0}+bz_1z_3^{h-w_1}+cz_2z_3^{h-w_2})-\varepsilon z_3^h$ and $\overline{\omega_{\nu,u,\varepsilon}}$ is independent of the choice of $l$ as long as $\partial F_{u,\varepsilon}/\partial z_l\ne0$ for any $l=1,2,3$. 
\end{proof}

\begin{proof}[Proof of Proposition \ref{prop:main}]
    The idea of the proof is to use the framing to construct a homotopy between the $\Delta(S_1,\ldots,S_\mu)$ and $\rho_{wt}$.

    Note that $\rho_{wt}$ is the path defined over the disc boundary corresponding to $\Tilde{f_0}$. The first step is to homotope $\rho_{wt}$ to a path defined over the disc boundary corresponding to $\Tilde{f}_1$. Similar to what has done in Example \ref{exp:weightedloop}, for each $\Tilde{f}_u$, the Gauss-Manin parallel transport identifies $H^2(M_{u,\mathrm{e}^{ih\theta}};\R)$ with $H^2(M_{u,1};\R)$. Therefore, the frame of the maximal positive definite space constructed in the previous lemma along the weighted path induces a path $B_u$ in $\mathcal{E}_u$, where $\mathcal{E}_u$ is the space of linear embeddings $e:\R^{b^+}\hookrightarrow H^2(M_{u,1};\R)$ whose image is a positive subspace. Since $\mathcal{E}_u$ and $\mathrm{SO}(b^+)$ are homotopy equivalent, then $[B_u]$ is indeed an element in $\pi_1(\mathrm{SO}(b^+))$ for each $u\in[0,1]$. By definition, $[B_0]=[\rho_{wt}]$. It is clear that $[B_1]=[\rho_{wt}]$.

    Let $B'_1$ be the loop of the concatenating of $A_i$, where $A_i$ is the loop going along the vanishing $\gamma_i$, followed by a small meridian loop $m_i$ around the critical value $s_i$ of $\Tilde{f}_1$, and going back along $\gamma_i$ for $i=1,\ldots,\mu$ (Figure \ref{fig:baseloop}). It is clear that the disc boundary circle is homotopic to $B_1'$. Then, using the frame for the maximal positive definite subspace constructed in the last lemma and the Gauss-Manin parallel transport, we regard $B_1'$ as a loop in $\mathcal{E}_1$. As a result, $B_1$ and $B_1'$ are homotopic in $\mathcal{E}_1$.

    \begin{figure}[H] 
    \centering 
    \includegraphics[width=0.4\textwidth]{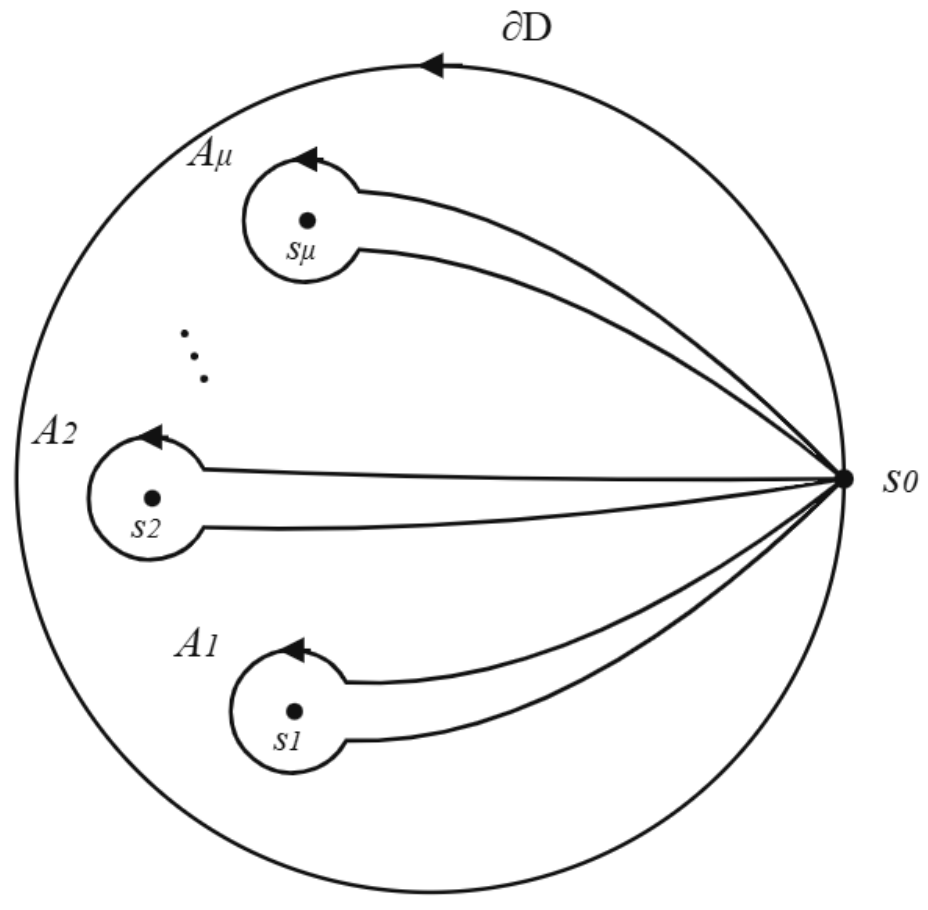}
    \caption{The loop concatenation.}\label{fig:baseloop}
    \end{figure}

    Recall that when defining $\Delta(S_1,\ldots,S_\mu)$, this path is first defined as a path of the framing of the maximal positive definite subspaces of $H^2(M;\R)$. Therefore, since $\mathcal{E}_1$ is homotopy equivalent to $SO(b^+)$, the proof of this proposition is complete if $\Delta(S_1,\ldots,S_\mu)$ and $B_1'$ are homotopic as paths in $\mathcal{E}_1$.

    To achieve this goal, we first show that $A_i$ is a path connecting $e_i$ and $\tau^*_{S_i}e_{i-1}=e_{i}$. By definition, $A_i$ is induced by transporting the elements from $H^2(M_{1,\varepsilon})$ to $H^2(M_{1,1})$ along the path defining $A_i$, where $\varepsilon$ is a point in the path. For $\alpha\in H^2(M_{1,1})$, the resulting element after the transporting is
    \begin{align*}
        P^{GM}_{\overline{\lambda_i}}(P^{GM}_{{m_i}}(P^{GM}_{{\lambda_i}}\alpha))=&P^{GM}_{\overline{\lambda_i}}(P^{GM}_{{\lambda_i}}\alpha+(P^{GM}_{{\lambda_i}}\alpha)(S'_i)\mathrm{PD}(S_i'))\\
        =&P^{GM}_{\overline{\lambda_i}}(P^{GM}_{{\lambda_i}}\alpha+\alpha(S_i)\mathrm{PD}(S_i'))\\
        =&\alpha+\alpha(S_i)\mathrm{PD}(S_i)\\
        =&\tau^*_{S_i}\alpha,
    \end{align*}
    where $S'_i$ is the vanishing sphere defining the vanishing cycle, $P^{GM}_\cdot$ stands for the Gauss-Manin parallel transport along the subscript path, and the second and the third equalities follows from the fact that the Gauss-Manin parallel transport preserved the pairing in $H^2$. After replacing $\alpha$ by the elements forming the frame $e_{i-1}$, it follows that $A_i$ is a path connecting $e_i$ and $\tau^*_{S_i}e_{i-1}=e_{i}$.

    Now $A_i$ and $\rho_t^i\circ e_{i-1}$ share the same ends. To show $\Delta(S_1,\ldots,S_\mu)$ are homotopic, one only needs to prove that the concatenating of $A_i$ and $\overline{\rho_t^i\circ e_{i-1}}$ (Proposition \ref{prop:deltloop}) is trivial as an element of $\pi_1\mathcal{E}_1$. Here the overline is refereed to reversing the direction of the path. Denote this loop by $C_i$. According to \cite[Lemma 2.6]{konno2025constraints}, $\mathcal{E}_1$ and $\mathcal{E}_1^i$ are homotopy equivalent using the family of isomorphisms $\rho_t^i$ with $i\in[0,\frac12]$, where $\mathcal{E}_1^i\subset\mathcal{E}_1$ consisting of embeddings of maximal positive definite subspaces whose image is orthogonal to $S_i$. Observe that $\rho_{\frac12}^i$ is the projection to the subspace orthogonal to $[S_i]$. Therefore, $[C_i]\in\pi_1(\mathcal{E}_1)$ if and only if $[\rho_\frac12^i(C_i)]\in\pi_1(\mathcal{E}_1^i)$ is trivial. A direct calculation shows that 
    \[
    \rho_\frac12^i\rho_t^i=\rho_\frac12^i,\ for\ all\ t,
    \]
    This implies that 
    \[
    \rho_\frac12^i(\rho_t^i\circ e_{i-1})=\rho_\frac12^i\circ  e_{i-1}
    \]
    is a constant path. As a result, $\rho_\frac12^i(A_i)$ is a loop and the next lemma shows that it is a trivial loop. So the proposition is proved since $\rho_\frac12^i(C_i)$ is the concatenating of $\rho_\frac12^i(\rho_t^i\circ e_{i-1})$ and $\rho_\frac12^i(A_i)$.

    The mod 2 formula is a direct consequence of Example \ref{exp:weightedloop} and Lemma \ref{lem:posframe}.
\end{proof}

\begin{lem}
    The loop $\rho_\frac12^i(A_i)$ is contractible.
\end{lem}

\begin{proof}
    Let $B\subset \C^3$ be a small ball around the critical point at the singular fiber such that $\Tilde{f}_1(B)\subset A_i$. For any regular value $\varepsilon$ near the critical value, define
    \[
    N_\varepsilon=M_{1,\varepsilon}-B,\ K_\varepsilon=\overline{M_{1,\varepsilon}-N_\varepsilon}.
    \]
    Since $\Tilde{f}_1$ is Morse, assume that the local coordinates near the critical points is chosen such that
    \[
    \Tilde{f}_1=z_0^2+z_1^2+z_2^2.
    \]
    Then $K_\varepsilon$ is the disc bundle of $T^*S^2$ \cite{ebeling2020distinguished}, and $\partial K_\varepsilon=\mathbb{RP}^3$. Here the 2-sphere in this $T^*S^2$ is the vanishing sphere \cite{ebeling2020distinguished}. By the Mayer–Vietoris sequence, for regular values $\varepsilon$
    \[
    H^2(M_{1,\varepsilon};\R)\cong H^2(K_\varepsilon;\R)\oplus H^2(N_\varepsilon;\R).
    \]
    Since the generator of $H^2(K_\varepsilon;\R)$ is generated by the vanishing sphere, $H^2(N_\varepsilon;\R)$ is the orthogonal completement of the vanishing cycle with respect to the intersection form. As a result, $\rho_\frac12^i([\mathrm{Im}\omega_{\nu,1,\varepsilon}])$ and $\rho_\frac12^i([\mathrm{Re}\omega_{\nu,1,\varepsilon}])$ can be regarded as elements in $H^2(N_\varepsilon;\R)$. Note that $N_\varepsilon$ is smooth when $\varepsilon$ is the critical value $s_i$, $\rho_\frac12^i([\mathrm{Im}\omega_{\nu,1,\varepsilon}])$ and $\rho_\frac12^i([\mathrm{Re}\omega_{\nu,1,\varepsilon}])$ can be extended to $H^2(N_{s_i};\R)$. In other words,
    \[
    \rho_\frac12^i([\mathrm{Im}\omega_{\nu,1,s_i}|_{N_{s_i}}]), \rho_\frac12^i([\mathrm{Re}\omega_{\nu,1,s_i}|_{N_{s_i}}])\in H^2(N_{s_i};\R),
    \]
    and are well-defined. Once we show that the space generated by the above class are positive definite, then $\rho^i_\frac12(A_i)$ is contractible.

    The strategy is prove that the space generated by $[\mathrm{Im}\omega_{\nu,1,s_i}|_{N_{s_i}}]$ and $[\mathrm{Re}\omega_{\nu,1,s_i}|_{N_{s_i}}]$ is positive definite and invariant under $\rho^i_\frac{1}{2}$. To show the positive definiteness, the argument is similar to the proof for Lemma \ref{lem:posframe}. First, $\omega_{\nu,1,s_i}$ can be extended to $Z_{1,s_i}$ across the resolution of the infinity divisor, where $Z_{1,s_i}$ is the resolution of the compactification of $M_{1,s_i}$ in the weighted projective space $\mathbb{CP}(w_0,w_1,w_2,1)$ at the infinity divisor by Lemma \ref{lem:extinfty}. The critical point $s_i$ is an $A_1$ singularity since locally the function can be written as 
    \begin{equation}\label{equ:coord}
        \Tilde{f}_1=z_0^2+z_1^2+z_2^2.
    \end{equation}
    Therefore, one time blow up at $s_i$ resolves this singularity. To show that $\omega_{\nu,1,s_i}$ can be extended again to get a holomorphic 2-form on $Z_{1,s_i}\#\overline{\mathbb{CP}^2}$, we perform the local calculations around $s_i$. By changing the coordinates again, assume $\Tilde{f}_1$ is expressed as 
    \[
    \Tilde{f}_1=xy-z^2.
    \]
    For calculating the blow up of the 2-forms, we only describe the calculation in the coordinate where $x\ne0$, the rest are similar. Let $(t,v,w)$ be the coordinates of the blow up at $s_i$, then
    \[
    x=t,y=uv,z=uw.
    \]
    So, $\Tilde{f}_1=t^2v-t^2w$. On the smooth fiber of $\Tilde{f}_1$, $v=w^2$. Therefore,
    \begin{align*}
        \omega_{\nu,1,s_i}=&\varphi_\nu(x,y,z)\frac{\dd x\wedge \dd z}{x}|_{\Tilde{f}_1=0}\\
        =&\varphi_\nu(t,tw^2,tw)\frac{t\dd t\wedge \dd w}{u}|_{\Tilde{f}_1=0}\\
        =&\varphi_\nu(t,tw^2,tw)\dd t\wedge \dd w|_{\Tilde{f}_1=0},
    \end{align*}
    where $\varphi_\nu$ is a holomorphic function induced by the coordinate change. It is clear that $\omega_{\nu,1,s_i}$ can be extended across the exceptional divisor holomorphically, denoted as ${\Omega}_{\nu,1,s_i}$. Because ${\Omega}_{\nu,1,s_i}$ generates a subspace of $H^{2,0}$ over a closed K\"ahler manifold $Z_{1,s_i}$, the same argument as in Lemma \ref{lem:posframe} implies that the space generated by 
    \[
    \mathrm{Re{\Omega}_{\nu,1,s_i}}|_{N_{s_i}}=\mathrm{Re{\omega}_{\nu,1,s_i}},\ \mathrm{Im{\Omega}_{\nu,1,s_i}}|_{N_{s_i}}=\mathrm{Im{\omega}_{\nu,1,s_i}}
    \]
    is positive definite.

    This lemma is proved if we show that
    \begin{equation}\label{equ:lim}
        \lim_{\varepsilon\rightarrow s_i}[\omega_{\nu,1,\varepsilon}]=[{\omega}_{\nu,1,s_i}],\ \forall\nu.
    \end{equation}
    In fact,
    \[
    \rho_\frac{1}{2}^i[\omega_{\nu,1,\varepsilon}]-[\omega_{\nu,1,\varepsilon}]=\frac12(\int_{S_i(\varepsilon)}\omega_{\nu,1,\varepsilon})\mathrm{PD}(S_i(\varepsilon)),
    \]
    where $S_i(\varepsilon)$ is defined as in Equation \ref{equ:vansphere} with $\gamma_i(t)$ replaced by $\varepsilon$ in the formula. It suffices to show that the above integral tends to 0 as $\varepsilon$ tends to $s_i$. In the coordinate of Equation \ref{equ:coord}, $s_i=0$ and
    \begin{align*}
        \int_{S_i(\varepsilon)}\omega_{\nu,1,\varepsilon}&=\int_{S_i(\varepsilon)}\frac{\varphi_\nu \dd z_1\wedge \dd z_2}{2z_0}\\
        &=\sqrt{|\varepsilon|}\int_{S_i(1)}\frac{\varphi_\nu \dd z_1\wedge \dd z_2}{2z_0},
    \end{align*}
    where function in the integral is written in the coordinate such that $z_0\ne0$ using the Poincar\'e residue map (hence is holomorphic), $\varphi_\nu$ is a holomorphic function induced by the local coordinates, and the second equation follows from a rescaling the integral sphere to make it unit radius. Therefore, the last integral above is finite and the whole equation vanishes as $\varepsilon\rightarrow0$. As a result, Equation \ref{equ:lim} is true and the proof of the lemma is complete.
\end{proof}

\section{Closed complex surfaces containing Brieskorn-Pham surfaces}\label{sec:comemb}

In this section, we construct closed complex surfaces with Brieskorn-Pham surfaces as submanifolds and then calculate their topological invariants.\\

Throughout the section, $X_\varepsilon$ standards for the space
\begin{equation*}
    \{[x:y:z_0:z_1:z_2]|\ \left\{
    \begin{aligned}
    \varepsilon x^b=&x^{b-a_0}z_0^{a_0}+x^{b-a_1}z_1^{a_1}+x^{b-a_2}z_2^{a_2}+z_0^b+z_1^b+z_2^b\\
    &+\beta x^{b-1}(y-x)\\
    y^c=&x^c+z_0^c+2z_1^c
    \end{aligned}
    \right.\}\subset\mathbb{CP}^4,
\end{equation*}
where $\beta$ is a transcendental number and $b,c>\max(a_0,a_1,a_2)$ are some fixed integers. The requirement of $\beta$ being transcendental guarantees that only $X_0$ has a single singular point for all small $|\varepsilon|$ and any integer tuple $(a_0,a_1,a_2,b,c)$. This is used in the next section dealing with smooth fillings. With a lengthy but easy calculation, one can show that $X_\varepsilon$ is smooth with small enough $|\varepsilon|\ne0$.

\begin{prop}\label{prop:MembX}
    Let $M_\varepsilon$ be a Brieskorm-Pham surface defined in Definition \ref{def:bpsurface}. Then $M_\varepsilon$ can be embedded into $X_\varepsilon$ with small $|\varepsilon|\ne0$. It is clear that $\tau_{S_i}$ can be extended trivially to $X_\varepsilon$ for $i=1,\ldots,\mu$.
\end{prop}

\begin{proof}
    Let $U$ be a small neighborhood of $(z_0,z_1,z_2)=(0,0,0)$. By choosing the branch such that $\rho(0,0)=1$, for $U$ small enough, the following is well defined:
    \[
    y:=\rho(z_0,z_1)=(1+z_0^c+2z_1^c)^\frac{1}{c}.
    \] 
    For $x=1$, the first equation defining $X_\varepsilon$ becomes
    \[
    z_0^{a_0}+z_1^{a_1}+z_2^{a_2}+z_0^b+z_1^b+z_2^b+\beta(y-1)=\varepsilon.
    \]
    The left hand side of the above equation is denoted as $g(z_0,z_1,z_2)$. Set $u=z_0^c+2z_1^c$ and define the holomorphic function
    \begin{equation*}
        q(u)=\left\{
        \begin{aligned}
            ((1+u)^{1/c}-1)/u,\ \ &u\ne0\\
            1/c,\qquad\qquad\qquad\quad &u=0
        \end{aligned}
        \right..
    \end{equation*}    
    The holomorphicity at zero follows by the Taylor expansion of $(1+u)^{1/c}$. Simple calculation shows
    \[
    \rho-1=uq(u).
    \]
    Consider the map $\Phi(z_0,z_1,z_2)=(\xi_0,\xi_1,\xi_2)$ defined by
    \begin{equation*}
        \left\{
        \begin{aligned}
            \xi_0=&z_0(1+z_0^{b-a_0}+\beta z_0^{c-a_0}q(u))^\frac{1}{a_0}\\
            \xi_1=&z_1(1+z_1^{b-a_1}+2\beta z_1^{c-a_1}q(u))^\frac{1}{a_1}\\
            \xi_2=&z_2(1+z_2^{b-a_2})^\frac{1}{a_2}
        \end{aligned}
        \right..       
    \end{equation*}
    Then, $\Phi:U\rightarrow V$ is biholomorphic for small neighborhoods $V$ and $U$ around $(0,0,0)$, since
    \begin{equation*}
        \dd\Phi(0,0,0)=\left[
        \begin{aligned}
        1\ \ 0\ \ 0\\
        0\ \ 1\ \ 0\\
        0\ \ 0\ \ 1
        \end{aligned}
        \right].
    \end{equation*}
    Let $f(\xi_0,\xi_1,\xi_2)=\xi_0^{a_0}+\xi_1^{a_1}+\xi_2^{a_2}$, so $g(z_0,z_1,z_2)=f(\Phi(z_0,z_1,z_2))$ by plugging the formulas. This shows that $\{g=\varepsilon\}\cap U$ and $\{f=\varepsilon\}\cap V$ are biholomorphic. For $x=1$, the second equation that defines $X_\varepsilon$ implies $y^c=1+z_0^c+2z_1^c$.
    By the above calculations, the embedding is given in the formula
    \begin{align*}
        \iota:M_\varepsilon&\rightarrow X_\varepsilon\\
        (\xi_0,\xi_1,\xi_2)&\mapsto[1:\rho(z_0,z_1):z_0:z_1:z_2],
    \end{align*}
    where $z_0,z_1,z_2$ are functions of $\xi_0,\xi_1,\xi_2$ using $\Phi$. It is clear that $\iota$ is injective, so it gives an embedding of $M_\varepsilon$ into $X_\varepsilon$.
\end{proof}

To apply Corollary 1.3 of \cite{konno2025constraints}, we compute the topological invariants of $X_\varepsilon$.

\begin{prop}\label{prop:topinf}
    Let $h\in H^2(\mathbb{CP}^4)$ be the generator and $H:=h|_{X_\varepsilon}\in H^2(X_\varepsilon)$.Then the complex manifold $X_\varepsilon$ is simply connected with
    \begin{itemize}
        \item $\sigma(X_\varepsilon)=\frac{bc}{3}(5-b^2-c^2)$,
        \item $c_1(\mathfrak{s}_{can})=(5-b-c)H$, where $\mathfrak{s}_{can}$ is the canonical $spin^c$ stucture induced by the complex structure on $TX_\varepsilon$,
        \item $c_1(\mathfrak{s}_{can})|_{M_\varepsilon}=0$,
        \item for $b=5+32r>5$ and $c=96l>0$, then $b^+>1$, $32|\sigma(X_\varepsilon)$, $32|c_1(\mathfrak{s}_{can})$, the Seiberg-Witten expected dimension $d(\mathfrak{s}_{can})=0$, and $SW(X_\varepsilon,\mathfrak{s}_{can})=\pm1$. Moreover, $X$ is a minimal surface of general type.
    \end{itemize}
\end{prop}

\begin{proof}
    By applying the Veronese embedding and the Lefschetz hyperplane theorem, we have $\pi_1(X_\varepsilon)=1$. Therefore, $b_1(X_\varepsilon)=b_3(X_\varepsilon)=0$. Recall that $b_0(X_\varepsilon)=b_4(X_\varepsilon)=1$, the Euler number is determined by $b_2(X_\varepsilon)$:
    \[
    \chi(X_\varepsilon)=2+b_2(X_\varepsilon).
    \]
    We first calculate the Chern class of $X_\varepsilon$ Since $X_\varepsilon$ is the intersection of degree $b$ and $c$ equations, Lemma 9.3.3 \cite{Waldron2020} implies that the normal bundle $N_{X_\varepsilon}$ is isomorphic to $\mathcal{O}(b)|_{X_\varepsilon}\oplus\mathcal{O}(c)|_{X_\varepsilon}$, where $\mathcal{O}(b)$ and $\mathcal{O}(c)$ are the holomorphic line bundle over $\mathbb{CP}^4$ of degree $b$ and $c$, respectively. As a result,
    \begin{equation}\label{equ:chernnormal}
        c(N_{X_\varepsilon})=(1+bH)(1+cH).
    \end{equation}
    By the exact sequence
    \[
    0\rightarrow TX_\varepsilon\rightarrow T\mathbb{CP}^4|_{X_\varepsilon}\rightarrow N_{X_\varepsilon}\rightarrow0,
    \]
    the Chern classes of the terms in the sequence satisfies
    \[
    c(TX_\varepsilon)c(N_{N_\varepsilon})=c(T\mathbb{CP}^4|_{X_\varepsilon})=(1+h)^5|_{X_\varepsilon}=(1+H)^5.
    \]
    Combining this with Equation \ref{equ:chernnormal}, we have
    \begin{align*}
        c(TX_\varepsilon)=&\frac{(1+H)^5}{(1+bH)(1+cH)}\\
        =&(1+H)^5(1-bH+b^2H^2)(1-cH+cH^2)\\
        =&1+(5-b-c)H+(b^2+bc+c^2-5b-5c+10)H^2.
    \end{align*}
    That is to say,
    \[
    c_1(TX_\varepsilon)=(5-b-c)H,\ c_2(TX_\varepsilon)=(b^2+bc+c^2-5b-5c+10)H.
    \]
    Note that the first Chern class of the canonical bundle
    \[
    c_1(K_{X_\varepsilon})=-c_1(TX_\varepsilon)=(b+c-5)H.
    \]
    Consequently, whenever $b+c>5$, $K_{X_\varepsilon}$ is ample and $X_\varepsilon$ is a minimal surface of general type. 
    
    Since $[X_\varepsilon]=bch^2$ in $H^2(\mathbb{CP}^4)$, the action of $H^2$ on $[X_\varepsilon]$ is $bc$. Therefore,
    \[
    b_2(X_\varepsilon)=\chi(X_\varepsilon)-2=c_2(TX_\varepsilon)([X_\varepsilon])-2=bc(b^2+bc+c^2-5b-5c+10)-2.
    \]
    Moreover, the signature
    \[
    \sigma(X_\varepsilon)=\frac13(c_1(TX_\varepsilon)^2-2c_2(TX_\varepsilon))[X_\varepsilon]=\frac{bc}{3}(5-b^2-c^2).
    \]
    Then $b_2(X_\varepsilon)$ and $\sigma(X_\varepsilon)$ implies
    \[
    b^+(X_\varepsilon)=\frac{bc}{6}(2b^2+3bc+2c^2-15b-15c+35)-1.
    \]
    For $c_1(\mathbf{\mathfrak{s}}_{can})$,
    \[
    c_1(\mathbf{\mathfrak{s}}_{can})=c_1(TX_\varepsilon)=(5-b-c)H.
    \]
    Note that $H$ can be represented by the intersection of the hyperplane $x=0$ and $X_\varepsilon$. Because this intersection is disjoint with $M_\varepsilon$,
    \[
    c_1(\mathbf{\mathfrak{s}}_{can})|_{M_\varepsilon}=0.
    \]
    Finally, the last point is obtained by direct calculations and applying Taubes' theorem \cite{Taubes1994} for the Seiberg-Witten invariants of symplectic manifolds.
\end{proof}

The main theorem is a consequence of Theorem B of \cite{konno2025constraints}, stated as below:

\begin{them}[{\cite[Theorem B]{konno2025constraints}}]\label{thm:kllmme}
Let $(X,\mathfrak{s})$ be a closed simply-connected
spin-$c$ smooth $4$-manifold. Let $S_1,\ldots,S_n$
be a collection of smoothly embedded $(-2)$-spheres.
Assume the following conditions hold:
\begin{itemize}
    \item both $c_1(\mathfrak{s})$ and $\sigma(X)$
    are divisible by $32$,
    \item $d(\mathfrak{s})=0$ and the Seiberg--Witten invariant $\operatorname{SW}(X,\mathfrak{s})$ is odd,
    \item $S_i$ pairs trivially with $c_1(\mathfrak{s})$, i.e., $c_1(\mathfrak{s})\cdot S_i=0$,
    \item The composition $\tau_{S_1}\cdots\tau_{S_n}$
    is smoothly isotopic to the identity.
\end{itemize}
Then $\Delta(S_1,\ldots,S_n)\equiv0\pmod{2}.$
\end{them}

\begin{them}\label{thm:exodiff0.5}
    For complex manifold $X_\varepsilon$ with integers $b=5+32r>5$ and $c=96m>0$, if 
    \[
    \Sigma_{\nu\in\mathcal{I}}N_\nu\equiv1\pmod{2},
    \]
    then $(\tau_{S_1}\cdots\tau_{S_\mu})^h$ is an exotic diffeomorphism of $X_\varepsilon$. Here $S_i$ are the vanishing spheres of $M_\varepsilon$,
    \[
    \mathcal{I}:=\{(i,j,k)|\, 1\le i\le a_0, 1\le j\le a_1,1\le k\le a_2, \frac{i}{a_0}+\frac{j}{a_1}+\frac{k}{a_2}<1\},
    \]
    and $N_\nu=iw_0+jw_1+kw_2-h$.
\end{them}

\begin{proof}
    By Proposition \ref{prop:topinf}, $X_\varepsilon$ satisfies all the topological requirements of \cite[Theorem B]{konno2025constraints}. However, $\Sigma_{\nu\in\mathcal{I}}N_\nu\equiv1\pmod{2}$. This shows that $(\tau_{S_1}\cdots\tau_{S_\mu})^h$ cannot be smoothly isotopic to the identity.

    It remains to show that $(\tau_{S_1}\cdots\tau_{S_\mu})^h$ act trivially on $H^2(X_\varepsilon;\Z)$. Once this is established, \cite{quinn1986isotopy} implies that $(\tau_{S_1}\cdots\tau_{S_\mu})^h$ is topologically isotopic to the identity. Recall that in the definition of Brieskorn-Pham surfaces, $a_0,a_1,a_2$ are assumed to be coprime to each other. Therefore, $\partial M_\varepsilon$ is a $\Z$-homology sphere by \cite[Proposition 2.2]{martin2023normal}. According to the Mayer–Vietoris sequence,
    \[
    H^2(X_\varepsilon;Z)\cong H^2(M_\varepsilon;\Z)\oplus H^2(\overline{X_\varepsilon-M_\varepsilon};\Z).
    \]
    As a result, $(\tau_{S_1}\cdots\tau_{S_\mu})^h$ acts trivially on $H^2(X_\varepsilon)$ since $(\tau_{S_1}\cdots\tau_{S_\mu})^h$ already induces a trivial action on $H^2(M_\varepsilon;\Z)$.
\end{proof}

\noindent Table \ref{tab:exotic-examples} are examples calculated for different tuples $(a_0,a_1,a_2,b,c)$ determining $X_\varepsilon$.

\begin{table}[htbp]
    \centering
    \renewcommand{\arraystretch}{1.3}
    \setlength{\tabcolsep}{15pt}
    \begin{tabular}{|c|c|c|}
        \hline
        $(a_0,a_1,a_2,b,c)$
        & $\displaystyle\sum_{\nu\in I}N_\nu$
        & Exotic \\
        \hline
        $(2,3,11,37,96)$   & $-5$     & Yes \\
        \hline
        $(2,3,13,37,96)$   & $-8$     & Undetermined \\
        \hline
        $(2,3,19,37,96)$   & $-21$    & Yes \\
        \hline
        $(2,3,23,37,96)$   & $-33$    & Yes \\
        \hline
        $(2,7,11,37,96)$   & $-105$   & Yes \\
        \hline
        $(2,7,15,37,96)$   & $-217$   & Yes \\
        \hline
        $(2,7,19,37,96)$   & $-369$   & Yes \\
        \hline
        $(2,7,23,37,96)$   & $-561$   & Yes \\
        \hline
        $(3,4,5,37,96)$    & $-14$    & Undetermined \\
        \hline
        $(3,7,10,37,96)$   & $-405$   & Yes \\
        \hline
        $(3,7,11,37,96)$   & $-505$   & Yes \\
        \hline
        $(3,7,19,37,96)$   & $-1701$  & Yes \\
        \hline
        $(3,10,11,37,96)$  & $-1185$  & Yes \\
        \hline
        $(3,11,14,37,96)$  & $-2535$  & Yes \\
        \hline
        $(6,7,11,37,96)$   & $-3475$  & Yes \\
        \hline
        $(7,10,11,37,96)$  & $-11565$ & Yes \\
        \hline
    \end{tabular}
    \caption{Examples of smooth projective complete
    intersections $X_\varepsilon$ and the exoticness
    of the diffeomorphisms $(\tau_{S_1}\circ\cdots\circ\tau_{S_\mu})^h$,
    where $h=\mathrm{lcm}(a_0,a_1,a_2)$. The associated Brieskorn--Pham singularities are not exceptional unimodal singularities.}
    \label{tab:exotic-examples}
\end{table}

\section{An obstruction to the smooth filling}\label{sec:obsfill}

The goal of this section is to show that the spin number is an obstruction of the existence of smooth fillings. Let $\mathcal{X}$ be the total space of the parametrized family $\{X_\varepsilon\}_{\varepsilon}$. It is prove in the last section that the projection $\pi^*:\mathcal{X}^*:=\mathcal{X}-X_0\rightarrow D^*$ that maps $X_\varepsilon$ to $\varepsilon$ is a fiber bundle with smooth fiber $X_\varepsilon$, where $D^*\subset\C$ is a disc with small radius punctured at $0$. We aim at proving the two theorems below.

\begin{them}\label{them:smX}
    The space $\mathcal{X}$ is smooth and the smooth fibration 
    \[
    \pi:\mathcal{X}\rightarrow D
    \]
    into the disc with small radius has only one singular fiber $X_0$. Moreover, there is a single singular point in $X_0$.
\end{them}

\begin{defn}\label{def:basechange}
    The $h$-order base change of $\mathcal{X}$ is the pullback fibration ${\{pi^{(h)}}^*:{\mathcal{X}^{(h)}}^*\rightarrow D^*$ over the punctured disc $D^*$
    \[
    \begin{tikzcd}
    {\mathcal{X}^{(h)}}^*
    \arrow[r, ]
    \arrow[d, "{\pi^{(h)}}^*"']
    &
    \mathcal{X}
    \arrow[d, "\pi"]\\
    D^* \arrow[r, "s^h"']
    &
    D^*
    \end{tikzcd}.
    \]
\end{defn}

\begin{them}\label{them:basechangefill}
    Let ${\pi^{(h)}}^*:{\mathcal{X}^{(h)}}^*\rightarrow D^*$ be $h$-order base change of $\pi:\mathcal{X}\rightarrow D$.
    Then, the fiber bundle
    \[
    \pi^{(h)}:{\mathcal{X}^{(h)}}^*\rightarrow D^*
    \]
    admits a topological filling. Moreover, if the spin number
    \[
    [\Delta(S_1,\ldots,S_\mu)^h]\equiv1\pmod2,
    \]
    then ${\mathcal{X}^{(h)}}^*$ has no smooth filling, here $S_i$ is the vanishing sphere of the Brieskorn-Pham surface $M_\varepsilon$ embedded into $X_\varepsilon$.
\end{them}

%\noindent The following technical lemma is needed for determining the singular point of $X_0$.
%\begin{lem}\label{lem:derivation}
%    Let $K$ be a finitely generated field extension of $\mathbb Q$ and let $\beta\in K$ be transcendental over $\mathbb Q$. Then there exists a derivation
%    \[
%    D:K\longrightarrow K
%    \]
%    such that $D(\beta)=1$ and $D|_{\mathbb Q}=0$.
%\end{lem}

%\begin{proof}
%    Extend $\{\beta\}$ to a transcendence basis $\{\beta,t_1,\ldots,t_s\}$ of $K$ over $\mathbb Q$. On the rational function field
%    \[
%    K_0=\mathbb Q(\beta,t_1,\ldots,t_s),
%    \]
%    define $D=\partial/\partial\beta$, and $D(t_1)=0$ for $i=1,\ldots,s$.

%    Since $K/K_0$ is a finite separable extension,
%   $D$ extends uniquely to $K$. Indeed, if $\alpha$
%is algebraic over $K_0$ with minimal polynomial
%$P(T)=\sum_{j=0}^m p_jT^j$, the extension is
%determined by
%\[
%D(\alpha)
%=-\frac{\sum_{j=0}^m D(p_j)\alpha^j}
%{P'(\alpha)},
%\]
%where $P'(\alpha)\neq0$ by separability.
%The resulting derivation satisfies
%$D(\beta)=1$ and $D|_{\mathbb Q}=0$.
%\end{proof}

\begin{proof}[Proof of Theorem \ref{them:smX}]
The smoothness of $\mathcal{X}$ and the smoothness of $X_0$ with $x=0$ are both a consequence of lengthy calculation of derivatives, so the proof is omitted here. 

It remains to show that there in only one singular point $[1:1:0:0:0]$ in $X_0$ when $x\ne0$. In the homogeneous coordinates, letting $x=1$, set
\[
g_j(z_j)=z_j^{a_j}+z_j^b,\qquad j=0,1,2.
\]
Then $X_0$ is defined by
\begin{equation}\label{eq:central-fiber}
\begin{cases}
\sum_{j=0}^2 g_j(z_j)+\beta(y-1)=0\\
y^c-1-z_0^c-2z_1^c=0
\end{cases}.
\end{equation}
The corresponding Jacobian matrix is
\[
J=
\begin{pmatrix}
\beta & g_0'(z_0) & g_1'(z_1) & g_2'(z_2)\\
cy^{c-1} & -cz_0^{c-1} & -2cz_1^{c-1} & 0
\end{pmatrix}.
\]
It is clear that $J$ is of rank one at $[1:1:0:0:0]$.

Suppose that $J$ is not full rank at $q=[1:y:z_0:z_1:z_2]\in X_0$. Since its second row cannot vanish on $X_0$ and $\beta\neq0$, we have $y\neq0$ and
\begin{equation}\label{eq:rank-deficiency}
\begin{cases}
g_0'(z_0)=-\dfrac{\beta z_0^{c-1}}{y^{c-1}}\\
g_1'(z_1)=-\dfrac{2\beta z_1^{c-1}}{y^{c-1}}\\
g_2'(z_2)=0
\end{cases}.
\end{equation}
%By Lemma \ref{lem:derivation}, the field $K=\mathbb Q(\beta,y,z_0,z_1,z_2)$ admits a derivation $D$ satisfying $D(\beta)=1$ and $D|_{\mathbb Q}=0$.
By \cite[Chapter 8, Corollary 5.6]{Lang02} and $\mathrm{ch}(\Q)=0$, there is a derivation $D$ such that $D(\beta)\ne0$ of the separable field extension $K:=\Q(\beta,y,z_0,z_1,z_2)$ over $\Q$. Then applying $D$ to Equation \ref{eq:central-fiber} gives
\begin{equation*}
\begin{cases}
\sum_{j=0}^2 g_j'(z_j)D(z_j)+\beta D(y)+D(\beta)(y-1)=0\\
y^{c-1}D(y)-z_0^{c-1}D(z_0)-2z_1^{c-1}D(z_1)=0
\end{cases}.
\end{equation*}
Plug Equation \ref{eq:rank-deficiency} into the first one above, then
\[
0=\frac{\beta}{y^{c-1}}(D(y)-z_0^{c-1}D(z_0)-2z_1^{c-1}D(z_1))+D(\beta)(y-1)=D(\beta)(y-1).
\]
Thus $y=1$, and Equation \ref{eq:central-fiber} can be reduced to
\begin{equation}\label{eq:y-one}
\begin{cases}
g_0(z_0)+g_1(z_1)+g_2(z_2)=0,\\
z_0^c+2z_1^c=0.
\end{cases}
\end{equation}
Suppose that $z_0,z_1\neq0$. Setting $\alpha=z_0/z_1$, we have $\alpha^c=-2$ and
$|\alpha|>1$. In particular, $\alpha$ is algebraic over $\mathbb Q$. Moreover, $g_2'(z_2)=0$ implies that $z_2$ is algebraic over $\mathbb Q$. Then replacing $z_0$ by $\alpha z_1$ into Equation \ref{eq:y-one}, it gives
\[
(\alpha^b+1)z_1^b+\alpha^{a_0}z_1^{a_0}+z_1^{a_1}+g_2(z_2)=0.
\]
Since $b>\max\{a_0,a_1\}$ and $|\alpha|>1$, this is a nonzero polynomial equation in $z_1$ with algebraic coefficients. Hence $z_1$ and $z_0=\alpha z_1$ are algebraic over $\mathbb Q$. However, Equation \ref{eq:rank-deficiency} implies
\[
\beta=-\frac{g_0'(z_0)}{z_0^{c-1}}
=-a_0z_0^{a_0-c}-bz_0^{b-c},
\]
contradicting the transcendence of $\beta$. Therefore, $z_0=z_1=0$. Finally, the remaining equation $g_2(z_2)=g_2'(z_2)=0$ implies $z_2=0$. If this is not true, then
\[
z_2^{b-a_2}=-1=-\frac{a_2}{b}.
\]
This contradicts to $b>a_2$. As a result, $q=[1:1:0:0:0]$. So $[1:1:0:0:0]$ must be the unique singular point of $X_0$.    
\end{proof}

\begin{prop}
Let
\[
\pi^*:\mathcal{X}^*\longrightarrow D^*
\]
be a proper smooth fiber bundle with closed oriented fiber $X$, where $D^*=D\setminus\{0\}$. As the fiber travels around $\partial D$ once, we get a self-diffeomorphism $\psi$ of the starting fiber. Then the following statements are equivalent:
\begin{enumerate}
    \item The family admits a differentiable smooth filling
    over $D$.
    \item The class
    \[
    [\psi]=1
    \quad\text{in }\pi_0\mathrm{Diff}(X).
    \]
\end{enumerate}
\end{prop}

\begin{proof}
$(1)\Rightarrow(2)$. Suppose that the family admits a smooth
filling
\[
\widetilde{\pi}:\widetilde{\mathcal{X}}\longrightarrow D.
\]
By Ehresmann's fibration theorem, $\widetilde{\pi}$ is a smooth fiber bundle. Since $D$ is contractible, this bundle is smoothly trivial. Restricting to $D^*$ gives
\[
\mathcal{X}^*\cong X\times D^*
\]
as smooth fiber bundles. Under this isomorphism, $\mathcal{X}^*$ can be smoothly extended to $D$. Since $D$ is contractible, $[\psi]=1$ in $\pi_0\mathrm{Diff}(X)$. \\

\noindent $(2)\Rightarrow(1)$. Since $D^*\cong S^1\times(0,1)$, the bundle $\mathcal{X}^*$ is smoothly isomorphic to
\[
T_\psi\times(0,1),
\]
where
\[
T_\psi=\frac{X\times[0,1]}{(x,1)\sim(\psi(x),0)}.
\]
Since $[\psi]=1$, there exists a smooth isotopy $\{\psi_s\}_{s\in[0,1]}$ from ${id}_X$ to $\psi$.
The map
\[
T_\psi\rightarrow X\times S^1,\qquad[x,s]\mapsto(\psi_s(x),\mathrm{e}^{2\pi\mathrm{i}s})
\]
is a well-defined smooth fiber-preserving diffeomorphism. Consequently,
\[
\mathcal{X}^*\cong T_\psi\times(0,1)\cong X\times D^*
\]
as smooth fiber bundles. It is clear that $X\times D^*$ can be extended across the origin. Therefore, $\mathcal{X}^*$ has a smooth filling.
\end{proof}

\begin{proof}[Proof of Theorem \ref{them:basechangefill}]
    The $h$-order base change space being a fiber bundle follows directly from Theorem \ref{them:smX}.
    
    To show that $\mathcal{X}^*$ admits a topological filling, note that the above proposition and its hold after replacing every word ``smooth" by ``topological."  By the definition of $\mathcal{X}^*$, the self-diffeomorphism of the fiber induced by traveling around $\partial D$ once is $(\tau_{S_1}\circ\cdots\circ\tau_{S_\mu})^h$. This is topologically isotopic to the identity according to Theorem \ref{thm:exodiff0.5}.

    Finally, combining Proposition \ref{prop:main} and Theorem \ref{thm:exodiff0.5}, $[(\tau_{S_1}\circ\cdots\circ\tau_{S_\mu})^h]\ne1$ in $\pi_0\mathrm{Diff}(X_\varepsilon)$ if $[\Delta(S_1,\ldots,S_\mu)^h]\equiv1\pmod2$. The above theorem implies that $\mathcal{X}^*$ admits no smooth fillings.
\end{proof}

\section{Symplectic embeddings and symplectic Torelli groups for surfaces of general type}\label{sec:sympemb}

The embeddings in the last section are not necessarily symplectic. The goal of this section is to modify the symplectic forms on $M_\varepsilon$ such that $M_\varepsilon\hookrightarrow X_\varepsilon$ is a symplectic where the target is equipped with the Fubini-Study form restricted form the ambient $\mathbb{CP}^4$.\\

\begin{prop}
    Let $M_\varepsilon$ be a Briskorn-Pham surface and $X_\varepsilon$ defined as 
    \begin{equation*}
    \{[x:y:z_0:z_1:z_2]|\ \left\{
    \begin{aligned}
    \varepsilon x^b=&x^{b-a_0}z_0^{a_0}+x^{b-a_1}z_1^{a_1}+x^{b-a_2}z_2^{a_2}+z_0^b+z_1^b+z_2^b\\
    &+\beta x(y-x)\\
    y^c=&x^c+z_0^c+2z_1^c
    \end{aligned}
    \right.\}\subset\mathbb{CP}^4,
\end{equation*}
    Then $M_\varepsilon$ can be equipped with a symplectic form $\omega_M$ that can be symplectic embedded into $X_\varepsilon$
    \[
    \hat\iota(M_\varepsilon,\lambda\omega_M)\rightarrow(X_\varepsilon,\omega_{FS}),\ \lambda>0
    \]
    where $\omega_M$ and $\omega_{FS}$ are the restrictions of the standard symplectic form from $\C^3$ and the Fubini-Study form from $\mathbb{CP}^4$, respectively. Moreover, $\hat\iota$ is smoothly isotopic to the embedding in Proposition \ref{prop:MembX}.
\end{prop}

\begin{proof}
    As in the construction of the embedding $\iota:M_\varepsilon\rightarrow X_\varepsilon$ in Proposition \ref{prop:MembX}, let $g(z_0,z_1,z_2)=z_0^{a_0}+z_1^{a_1}+z_2^{a_2}+z_0^{b}+z_1^{b}+z_2^{b}$. Moreover, there is a biholomorphic function $\Phi:U\rightarrow V$ such that $g=f\circ\Phi$, where $U,V$ are small neighborhood of $(0,0,0)\in\C^3$ and $f(\xi_0,\xi_1,\xi_2)=\xi_0^{a_0}+\xi_1^{a_1}+\xi_2^{a_2}$. In other words, $(g=\varepsilon)\cap V$ and $(f=\varepsilon)\cap U$ are biholomorphic through $\Phi$.

    We next modify the ambient symplecitic form so that this embedding becomes symplectic. Choose an open ball $B^3_R(0)$ whose closure is contained in $\Phi^{-1}(V)$. By \cite[Lemma 1.7]{seidel2003long}, there exist a symplectic form $\widetilde\omega$ on $B^3_R(0)$, a positive constant $\lambda$, and a smaller radius $r>0$, such that $\widetilde\omega$ is the restriction of the Fubini-Study form pulled back by $\Phi$ near $\partial B^3_R(0)$, and $\widetilde\omega=\lambda\omega_{\mathbb C^3}$ in $B^3_r(0)$, where $\omega_{\mathbb C^3}$ is the standard symplectic form on $\C^3$. One can choose $\delta$ small in Definition \ref{def:bpsurface} such that $M_\varepsilon\subset B^3_r(0)$. Then 
    \[
    \iota:(M_\varepsilon,\lambda\omega_M)\rightarrow(X_\varepsilon,\omega)
    \]
    is a symplectic embedding, where
    \[
     \omega=(\Phi^{-1})^*(\widetilde\omega)=
     \begin{cases}
        (\Phi^{-1})^*\omega_{\C^3}|_{X_\varepsilon},&\text{on } X_\varepsilon\cap B^3_R(0),\\
        \omega_{\mathrm{FS}},&\text{outside }B^3_R(0).
    \end{cases}
    \]
    Next we prove that $[\omega]=[\omega_{\mathrm{FS}}]$. Because $\omega-\omega_{\mathrm{FS}}$ is a compactly supported closed two-form on $B^3_R(0)$, and
    \[
        H_c^2(B^3_R(0);\mathbb R)=0,
    \]
    there exists a compactly supported one-form $\beta$ on $B^3_R(0)$ such that
    \[
        \omega-\omega_{\mathrm{FS}}=d\beta.
    \]
    Since $\beta$ can be extended trivially outside the ball, $[\omega]=[\omega_{\mathrm{FS}}]$.

    The final step is to apply the Moser's trick. Since both $\omega$ and $\omega_{\mathrm{FS}}$ are symplectic forms compatible with respect to the same complex structure, their convex combination
    \[
       \omega_t=(1-t)\omega_{\mathrm{FS}}+t\omega,\qquad t\in[0,1],
    \]
    is a family of symplectic forms representing the same cohomology class. By Moser's trick, there exists a diffeomorphism $\phi:X_\varepsilon\to X_\varepsilon$, smoothly isotopic to the identity, such that
    \[
       \phi^*\omega=\omega_{\mathrm{FS}}.
    \]
    As a result, $\hat\iota:=\phi^{-1}\circ\iota$ is the desired embedding.
\end{proof}

The proof of Theorem \ref{thm:mainexot} is at the hand.

\begin{proof}[Proof of Theorem \ref{thm:mainexot}]
    This is a direct consequence of the above proposition and Theorem \ref{thm:exodiff0.5}.
\end{proof}

\bibliographystyle{alpha}
\bibliography{references}

\noindent Zuyi Zhang, 
Beijing International Center for Mathematical Research, 
 Peking University, Beijing, 100871, China\\
\noindent Email: zhangzuyi1993@pku.edu.cn\\
\noindent Personal page: sites.google.com/view/zuyizhangmath/home

\end{document}